\documentclass[11pt]{article}

\usepackage[margin=1in]{geometry}
\usepackage{amsmath,amssymb,amsthm,mathtools}
\usepackage{booktabs}
\usepackage{graphicx}
\usepackage{xcolor}
\usepackage[numbers,sort&compress]{natbib}
\usepackage[hidelinks]{hyperref}
\usepackage{microtype}

\newtheorem{theorem}{Theorem}[section]
\newtheorem{proposition}[theorem]{Proposition}
\newtheorem{corollary}[theorem]{Corollary}
\newtheorem{lemma}[theorem]{Lemma}
\theoremstyle{definition}
\newtheorem{definition}[theorem]{Definition}
\newtheorem{remark}[theorem]{Remark}

\DeclareMathOperator{\Exp}{Exp}
\DeclareMathOperator{\Var}{Var}
\DeclareMathOperator{\KL}{KL}

\newcommand{\Pp}{\mathbb P}
\newcommand{\Ee}{\mathbb E}
\newcommand{\var}{\operatorname{Var}}
\newcommand{\R}{\mathbb R}

\newcommand{\pl}{\ell^{\mathrm p}}
\newcommand{\ip}[2]{\left\langle #1,#2\right\rangle}

\title{Nuisance structure and the detection boundary for population size:\\
how profiling can change a rate, not only a constant}
\author{Yakir Berchenko\\
\small Department of Industrial Engineering and Management\\
\small Ben-Gurion University of the Negev\\
\small \texttt{berchenk@bgu.ac.il}}
\date{August 2026}

\begin{document}
\maketitle

\begin{abstract}
A population of unknown size \(N\) consists of units that fail independently at an unknown common rate and are removed upon failure. In the Jelinski–Moranda model, inference is based on the first \(m\) inter-failure times, and the maximum likelihood estimate of \(N\) can be infinite even after many failures. We show that the critical number of failures required to distinguish a finite population from an inexhaustible one is of order \(N^{2/3}\). The exponent, however, is not a property of the model alone. Allowing an unknown linear trend in log testing effort absorbs the leading depletion signal and raises the critical scale to \(N^{4/5}\).

To explain this phenomenon, we develop a residual-information criterion for detectability in regular independent-observation experiments with expanding index designs and finite-dimensional nuisance trends. For parameter-scale perturbations of the form \(a(j/N)\), if the first Taylor term that cannot be absorbed by nuisance fitting has order \(q\) and the residual Fisher information for the index trend \(j^q\) grows as \(m^{2q+1}\), the critical scale is \(N^{2q/(2q+1)}\). Below this scale, suitably matched laws merge in total variation; above it, consistent separation is possible. At the critical scale, a one-sided Gaussian likelihood experiment emerges.

For the Jelinski–Moranda model, we provide the first asymptotic characterization of the probability of an infinite estimate: when \(m/N^{2/3}\to s\), this probability converges to \(\Phi(-s^{3/2}/\sqrt{12})\). We also derive the estimator's limiting distribution, a feasible locally optimal test and likelihood confidence sets that may be unbounded.
\end{abstract}



\section{Introduction}

In regular likelihood problems, estimating nuisance parameters typically reduces the information available for inference without changing the asymptotic rate, provided that the efficient information remains positive. This paper examines a different situation: one in which the departure from a limiting model is so closely confounded with nuisance variation that estimating an additional nuisance parameter changes the asymptotic order of the amount of data required for detection. The phenomenon arises naturally in a classical problem of population-size estimation.

Suppose that a system initially contains an unknown number $N$ of faults. Each fault independently generates failures at an unknown common rate $\phi>0$ and is removed upon its first detection. In the Jelinski-Moranda (JM) model of software reliability \citep{jelinski1972}, the observations are the first $m$ inter-failure times, which are independent and satisfy
\[
 X_j\sim\Exp\{\phi(N-j)\},\qquad j=0,\ldots,m-1.
\]
Both $N$ and $\phi$ are to be estimated from these observations. Although the failure rate decreases as faults are removed, its decline may be difficult to distinguish from random variation in the observed waiting times. Writing the rate as $\lambda(1-j/N)$, where $\lambda=\phi N$ is the initial system-wide failure rate, makes the difficulty particularly apparent: for large $N$ relative to $m$, the rates are nearly constant. The limiting model, obtained by setting $1/N=0$ while retaining an unknown $\lambda$, is an inexhaustible source of failures with constant rate. The maximum likelihood estimate of $N$ need not be finite; the profile likelihood may attain its maximum at this boundary.

Consider, for instance, a system that actually contains $10{,}000$ faults, of which $1{,}000$ have been detected and removed. Despite having observed a tenth of the entire population, the probability that the maximum likelihood estimate of its initial size is infinite is approximately $18\%$, according to the asymptotic formula derived in this paper. This behaviour illustrates a difficulty more fundamental than that of obtaining a satisfactory finite estimate: how much information do the observations contain about whether the population is finite in the first place?

The existence of finite maximum likelihood estimates, the geometry of the likelihood and alternative inferential procedures for the JM model have been studied extensively \citep{littlewood1981,joe1985,lindsay1987,raftery1988,osborne2000}. Existing asymptotic analyses \citep{vanpul1992,osborne2000} have dealt primarily with regimes in which the number of observed failures is of order $N$. Our interest is in an earlier stage of the experiment, where $m/N\to0$ and only a vanishing fraction of the population has been observed. We ask how many failures are required before the decline in failure rate becomes statistically distinguishable from a constant rate.

We show that the critical number is of order $N^{2/3}$. Below this scale, the finite-population experiment and a suitably matched constant-rate experiment merge in total variation. Above it, finiteness can be detected consistently. At the critical scale, the suitably rescaled profile log-likelihood ratio converges to that of a nondegenerate one-sided Gaussian shift experiment, yielding, among other consequences, the approximation
\[
 \Pp_N(\widehat N=\infty)\approx
 \Phi\left(-\frac{m^{3/2}}{\sqrt{12}\,N}\right)
\]
used in the preceding example.

A further question arises when the model is allowed to accommodate changes in testing effort. Suppose that the event rates are now given by
\[
 \lambda_{j,N}=\exp\{\eta+\gamma(j/m)\}(1-j/N),
 \qquad j=0,\ldots,m-1,
\]
where $\eta$ is the unknown initial log rate and $\gamma$ is an unknown drift coefficient. The latter permits an increasing or decreasing trend in effort over the course of the experiment. On the log-rate scale, the effect of depletion is
\[
 \log(1-j/N)=-\frac{j}{N}-\frac{j^2}{2N^2}
 +O\{(j/N)^3\}.
\]
The linear term cannot be distinguished from a change in $\gamma$. Consequently, evidence of depletion must come from curvature in the log-rate profile rather than from its slope. We show that this apparently modest change in the nuisance model raises the critical number of observations from order $N^{2/3}$ to order $N^{4/5}$.

The distinction between the two nuisance parameters is instructive. Relative to a known initial rate, estimating the common level reduces the leading information about depletion by a factor of four, but leaves the detection exponent unchanged. Allowing an unknown linear trend, on the other hand, removes the first-order depletion signal altogether. The confounding itself is elementary; its consequence for the amount of data required for detection is less immediate. Why do these particular exponents arise, and what determines whether fitting a nuisance parameter changes only an information constant or the detection rate itself?

To answer these questions, we consider a class of independent-observation experiments in which a regular scalar-parametric family is observed along an expanding index design, with a departure of the form $a(\rho j)$ on the parameter scale and a fixed finite-dimensional space of nuisance trends, where $\rho=1/N$ in the population-size problem. In the regime $m\rho\to0$, the observations reveal an increasingly small initial segment of the curve $a$. The terms in its local expansion that lie in the nuisance space can be absorbed by changes in the nuisance coefficients; the remaining terms carry information about the departure. We call the order $q$ of the first term that survives this adjustment the \emph{effective contact order}.

For the index designs considered here, under the stated regularity and residual-information conditions, the squared Fisher norm of the first surviving departure is of order $\rho^{2q}m^{2q+1}$. In population-size terms this is $m^{2q+1}/N^{2q}$. Intuitively, a departure of order $q$ has a typical magnitude of order $(m/N)^q$ over the observed segment, and its squared magnitude accumulates over $m$ observations. The resulting critical scale is
\[
 m\asymp N^{2q/(2q+1)}.
\]
The exponent thus depends jointly on the local form of the departure and the nuisance space. The sharp constant is determined by the Fisher-orthogonal projection of the surviving term off that space. Although the critical observation scale varies with $q$, the normalized profile likelihood has the same one-sided Gaussian form at each of these scales. The framework therefore separates the rate at which information emerges from the form of the likelihood experiment that emerges at that rate.

The starting point for this analysis is classical efficient-score and profile-likelihood theory \citep{choi1996,murphy2000}, together with the use of local alternatives to describe the effects of nuisance specification \citep{saikkonen1989}. Higher-order identification and nonstandard rates have also been studied in models with singular Fisher information \citep{rotnitzky2000}, in second-order identified models \citep{dovononhall2018}, and in finite mixtures \citep{chen1995,nguyen2013,honguyen2016,honguyen2019,heinrich2018}. Here the individual observation family remains regular, but the sequence of experiments becomes weakly informative about the departure as the observed fraction shrinks and nuisance trends absorb its leading terms. The question is not whether higher-order identification can occur, but how it translates into a sharp observation requirement in this setting, and what the resulting likelihood implies for inference. The emergence of a nondegenerate experiment at a critical scale also connects this problem to the literature on detection boundaries \citep{ingstersuslina2003,donohojin2004}, while the form of the limiting likelihood is familiar from classical boundary asymptotics \citep{chernoff1954,self1987,brazzale2024}.

The paper makes three main contributions.  First, it develops a residual-information criterion for the class of expanding-index experiments described above, establishing asymptotic merging, a nondegenerate local likelihood experiment and consistent separation in the respective regimes. Second, it quantifies the effect of nuisance fitting, distinguishing the loss of a constant fraction of information from the loss of entire orders of
the departure.  Explicit information constants for polynomial nuisance trends and the example of unknown effort drift illustrate the distinction.
These results also point to a practical consideration in the design of depletion experiments: recording a potential nuisance covariate, where this allows its contribution to enter as a known offset, can change the order of
the number of observations needed for detection.  Third, for the JM model it identifies the sparse-regime detection boundary and derives the associated
likelihood experiment, the limiting distribution of the population-size estimator (including its atom at infinity), a feasible locally optimal test, and likelihood confidence sets that can be unbounded when the data do not exclude the limiting model.

The remainder of the paper proceeds from the JM model to the more general setting and then returns to population-size inference. Section~\ref{sec:jmmotivation} provides a direct calculation of the JM boundary. Sections~\ref{sec:general} and~\ref{sec:tangent} develop the theory for an unknown common level and for general finite-dimensional nuisance trends, respectively. Section~\ref{sec:exp} gives an exact profile likelihood identity for exponential rate arrays, and Section~\ref{sec:jm} derives the inferential consequences for JM. Section~\ref{sec:simulation} examines the finite-sample accuracy of the critical approximations through simulation. The final section discusses the implications of the results and possible extensions.

\section{A motivating calculation: Jelinski-Moranda model}
\label{sec:jmmotivation}

\subsection{Exact profile likelihood}

In the fixed-count JM experiment, the first $m$ inter-failure times satisfy
\begin{equation}
 X_j\sim\Exp\{\phi(N-j)\},
 \qquad j=0,\ldots,m-1,
\label{eq:jm}
\end{equation}
independently, where $N\ge m$ is an unknown integer and $\phi>0$ is unknown.
Write
\[
 S_m=\sum_{j=0}^{m-1}X_j,\qquad
 T_m=\sum_{j=0}^{m-1}jX_j,\qquad
 W_m=\frac{T_m}{S_m},\qquad
 \bar j_m=\frac{m-1}{2}.
\]
After profiling $\phi$ and putting $\rho=1/N$, the log likelihood relative to
the infinite-reservoir limit $\rho=0$ is
\begin{equation}
 \pl_m(\rho)-\pl_m(0)
 =
 \sum_{j=0}^{m-1}\log(1-\rho j)
 -
 m\log(1-\rho W_m).
\label{eq:jmprofile}
\end{equation}

The following criterion is classic
\citep{littlewood1981,joe1985,lindsay1987}; we include a short proof because
it connects the exact likelihood to the limiting calculation.

\begin{proposition}[Infinite-MLE criterion]
\label{prop:existence}
Apart from the probability-zero tie case,
\begin{equation}
 \widehat N=\infty
 \quad\Longleftrightarrow\quad
 W_m\le\bar j_m.
\label{eq:existence}
\end{equation}
\end{proposition}

\begin{proof}
Differentiating \eqref{eq:jmprofile} gives
\[
 \frac{d}{d\rho}\pl_m(\rho)
 =
 m\frac{W_m}{1-\rho W_m}
 -
 \sum_{j=0}^{m-1}\frac{j}{1-\rho j}.
\]
For fixed admissible $\rho$, the function
$x\mapsto x/(1-\rho x)$ is increasing and convex.  If
$W_m\le\bar j_m$, monotonicity and Jensen's inequality imply that the
derivative is nonpositive for every $\rho$, so the maximum is attained at
$\rho=0$.  If $W_m>\bar j_m$, the derivative at zero is positive, so
$\rho=0$ cannot maximize the likelihood.  Unimodality (see, e.g.,Lemma \ref{lem:unimodal} in the appendix) then gives a finite maximizer.
\end{proof}

\subsection{Where the two-thirds exponent comes from}

Define
\begin{equation}
 V_m
 =
 \sum_{j=0}^{m-1}(j-\bar j_m)^2
 =
 \frac{m(m^2-1)}{12}
\label{eq:Vm}
\end{equation}
and the scale-invariant statistic
\begin{equation}
 \Delta_m
 =
 \frac{m}{S_m\sqrt{V_m}}
 \sum_{j=0}^{m-1}(j-\bar j_m)X_j
 =
 \frac{m(W_m-\bar j_m)}{\sqrt{V_m}}.
\label{eq:Deltajm}
\end{equation}
In particular, Proposition~\ref{prop:existence} says that
$\widehat N=\infty$ precisely when $\Delta_m\le0$.

A second-order expansion of \eqref{eq:jmprofile} around $\rho=0$ gives,
on the scale relevant below,
\begin{equation}
 \pl_m(\rho)-\pl_m(0)
 =
 \rho\sqrt{V_m}\,\Delta_m
 -
 \frac{\rho^2V_m}{2}
 +o_{\Pp}(1).
\label{eq:jmheuristic}
\end{equation}
Thus the experiment is organized by the information parameter
\begin{equation}
 h_{m,N}
 =
 \frac{\sqrt{V_m}}{N}
 \sim
 \frac{m^{3/2}}{\sqrt{12}\,N}.
\label{eq:hjm}
\end{equation}
If $h_{m,N}$ vanishes, the linear and quadratic likelihood contrasts vanish;
if it diverges, depletion becomes unmistakable.  The first nondegenerate
scale is therefore
\[
 m^{3/2}/N\asymp1,
 \qquad\text{or equivalently}\qquad
 m\asymp N^{2/3}.
\]
Sections~\ref{sec:general}-\ref{sec:jm} turn this calculation into a general
theorem and then return to JM for sharp distributional consequences.

\section{Analytic contact: the common-level case}
\label{sec:general}

\subsection{Finite models approaching a limiting model}

Let $\{P_\theta:\theta\in\Theta\subset\R\}$ be a dominated statistical
model.  The symbol $P_\theta$ denotes the probability law of one
observation when its scalar parameter is $\theta$.  For each $m$, we
consider independent observations whose parameters change with the event
index:
\begin{equation}
 \begin{split}
  Y_{j,m}\mid(\eta,\rho)&\sim P_{\theta_{j,m}},\\
  \theta_{j,m}&=\eta+a(\rho j),
  \qquad j=0,\ldots,m-1.
 \end{split}
\label{eq:analyticmodel}
\end{equation}
Here $j$ indexes the design position, $\eta$ is an unknown common level, and
$\rho\geq0$ controls departure from the limiting model.  The function $a$
describes how the one-observation parameter changes with the rescaled design
coordinate $x=\rho j$; we take $a(0)=0$ without loss of generality and require
$\eta+a(\rho j)\in\Theta$.  At the boundary $\rho=0$,
\[
    \theta_{0,m}=\cdots=\theta_{m-1,m}=\eta,
\]
so the model reduces to $m$ iid observations from $P_\eta$.  In
finite-population depletion applications, $\rho=N^{-1}$, so that
$x=j/N$ is the fraction of the initial population already depleted.
We refer to this homogeneous boundary experiment as the limiting model.

The JM model gives an exact translation of every symbol:
\begin{equation}
 P_\theta=\Exp(\text{rate}=e^\theta),\qquad
 \eta=\log(\phi N),\qquad
 \rho=N^{-1},\qquad
 a(x)=\log(1-x).
\label{eq:jmmapping}
\end{equation}
Indeed,
$e^{\eta+a(\rho j)}=\phi N(1-j/N)=\phi(N-j)$, recovering
\eqref{eq:jm}.

\begin{definition}[Analytic contact order]
\label{def:contact}
Suppose that $a$ is analytic near zero and, for some integer $q\ge1$,
\begin{equation}
 a(x)=-cx^q+O(x^{q+1}),\qquad c>0.
\label{eq:contact}
\end{equation}
Equivalently, the first $q-1$ derivatives of $a$ vanish at zero and
$a^{(q)}(0)=-c q!\ne0$.  We then say that the model has \emph{analytic contact of order \(q\)}
with the limiting model.
\end{definition}

Thus ``contact'' is simply order of tangency at the boundary
$\rho=0$.  For $q=1$ the finite model separates from the limiting model
linearly; for $q>1$ the depletion curve is flatter and its first visible
departure occurs later.  The JM profile
$a(x)=\log(1-x)=-x+O(x^2)$ has contact order one.  The first polynomial
trend determines the information boundary.  Profiling the unknown common
level centers this trend and changes the sharp information constant, but not
its order.  Section~\ref{sec:tangent} gives the more general definition in
which a multidimensional nuisance space can absorb lower-order trends entirely.

\subsection{Elementary smoothness assumptions}

Write $\ell_\theta(y)$ for the log-density of $P_\theta$,
$\dot\ell_\theta(y)$ for its score, and
\[
  I(\theta)=\Ee_\theta\{\dot\ell_\theta(Y)^2\}\in(0,\infty).
\]
Fix a compact set $K_0$ in the interior of $\Theta$ and an open neighborhood
$U\supset K_0$.  In the common-level results below, $K_0$ may be an
arbitrarily
small neighborhood of $\eta_0$; in Section~\ref{sec:tangent}, it contains the
entire baseline parameter curve.  To keep the argument self-contained, we use
the following familiar sufficient assumptions uniformly on $U$.
\begin{enumerate}
\item[(S1)] For almost every $y$, the map
$\theta\mapsto\ell_\theta(y)$ is three times continuously differentiable on
$U$, the support does not depend on $\theta$, and differentiation under the
integral sign is valid for the expectations below.  The information identities
\[
 \Ee_\theta\dot\ell_\theta(Y)=0,
 \qquad
 -\Ee_\theta\ddot\ell_\theta(Y)
 =
 \Ee_\theta\{\dot\ell_\theta(Y)^2\}
 =
 I(\theta)
\]
hold there, and $I$ is continuously differentiable and bounded away from zero
on $U$.
\item[(S2)] For some $\delta>0$, there is a function $M$ such that
\[
 \sup_{\theta\in U}
 \left\{
  |\dot\ell_\theta(y)|^{2+\delta}
  +
  |\ddot\ell_\theta(y)|^{2+\delta}
  +
  |\ell_\theta^{(3)}(y)|
 \right\}
 \le M(y),
\qquad
 \sup_{\theta\in U}\Ee_\theta M(Y)<\infty.
\]
\item[(S3)] The log likelihood is concave in $\theta$, and the homogeneous
model is identifiable: its expected log likelihood has the unique maximizer
$\eta_0$ and is bounded away from its maximum outside every neighborhood of
$\eta_0$.
\end{enumerate}
Concavity is stronger than necessary, but it avoids a separate technical
assumption about localization of the profile maximizer.  It is satisfied by
the exponential-rate model below and by standard regular one-parameter
exponential families.  More general smooth models require only that the
relevant local maximizer is root-$m$ consistent.

\subsection{Main Theorem I}
Define the centred index trend
\begin{equation}
 g_{j,q,m}=j^q-\frac1m\sum_{i=0}^{m-1}i^q,
 \qquad
 V_{q,m}=\sum_{j=0}^{m-1}g_{j,q,m}^2.
\label{eq:gdef}
\end{equation}
The order of $V_{q,m}$ follows from a Riemann sum.  Writing
$j^q=m^q(j/m)^q$,
\[
 \frac{V_{q,m}}{m^{2q+1}}
 =
 \frac1m\sum_{j<m}
 \Bigl\{\Bigl(\frac jm\Bigr)^{\!q}
        -\frac1m\sum_{i<m}\Bigl(\frac im\Bigr)^{\!q}\Bigr\}^2
\]
is the empirical variance of $(j/m)^q$ over the uniform grid
$\{0,1/m,\ldots,(m-1)/m\}$, and therefore converges to the variance of
$T^q$ for $T$ uniform on $(0,1)$:
\begin{equation}
 v_q=\Var(T^q)
 =\int_0^1t^{2q}dt-\Bigl(\int_0^1t^qdt\Bigr)^{2}
 =\frac1{2q+1}-\frac1{(q+1)^2}
 =\frac{q^2}{(2q+1)(q+1)^2}.
\label{eq:vdef}
\end{equation}
Because the summand is polynomial in $j/m$, the Riemann error is
$O(m^{-1})$, giving
\begin{equation}
 V_{q,m}=v_qm^{2q+1}\{1+O(m^{-1})\}.
\label{eq:Vasym}
\end{equation}
For $q=1$, $v_1=1/12$ and \eqref{eq:Vasym} reproduces the exact JM
identity $V_m=m(m^2-1)/12$ of \eqref{eq:Vm}.

Write
\[
 \bar a_m(\rho)=\frac1m\sum_{j=0}^{m-1}a(\rho j),
 \qquad
 \widetilde a_{j,m}(\rho)=a(\rho j)-\bar a_m(\rho).
\]
Replacing $\eta$ by $\eta+\bar a_m(\rho)$ shows exactly how the nuisance
parameter absorbs the average perturbation.  The centered residual satisfies
\[
 \widetilde a_{j,m}(\rho)
 =
 -c\rho^q g_{j,q,m}
 +O\{\rho^{q+1}(j^{q+1}+m^{q+1})\}.
\]
Its natural local magnitude is therefore
\begin{equation}
  h_m=h_m(\rho,\eta)
  :=c\rho^q\{I(\eta)V_{q,m}\}^{1/2}.
\label{eq:hdef}
\end{equation}

For a sample $Y^{(m)}=(Y_{0,m},\ldots,Y_{m-1,m})$, write
\[
 \pl_m(\rho)=
 \sup_{\eta\in\Theta}
 \sum_{j=0}^{m-1}\ell_{\eta+a(\rho j)}(Y_{j,m})
\]
for the profile log-likelihood.  The next theorem is (a specific simpler version of) our central result.

\begin{theorem}[Profile likelihood with an unknown common level]
\label{thm:lan}
Suppose (S1)-(S3) and \eqref{eq:contact} hold, and let $m\to\infty$.
For $h\ge0$, define
\[
 \rho_m(h)
 =
 \left\{\frac{h}{c\{I(\eta_0)V_{q,m}\}^{1/2}}\right\}^{1/q}.
\]
Under the homogeneous experiment and uniformly for $h$ in compact subsets of
$[0,\infty)$,
\begin{equation}
  \pl_m\{\rho_m(h)\}-\pl_m(0)
  =
  h\Delta_m-\frac{h^2}{2}+o_{\Pp}(1),
\label{eq:profilelan}
\end{equation}
where
\begin{equation}
 \Delta_m
 =
 -\frac{\sum_{j=0}^{m-1}
        g_{j,q,m}\dot\ell_{\eta_0}(Y_{j,m})}
        {\{I(\eta_0)V_{q,m}\}^{1/2}}
 \ \Longrightarrow\ \mathcal N(0,1).
\label{eq:central}
\end{equation}
For a fixed \(h_0\geq0\), consider the local alternative matched to the limiting model,
\begin{equation}
 Y_{j,m}\sim
 P_{\eta_0+\widetilde a_{j,m}\{\rho_m(h_0)\}},
 \qquad j=0,\ldots,m-1.
\label{eq:matchedalternative}
\end{equation}
Under this alternative,
\begin{equation}
  \Delta_m\Longrightarrow\mathcal N(h_0,1).
\label{eq:alternativeclt}
\end{equation}
Consequently the profiled likelihood is asymptotically that of a single
observation $\Delta\sim\mathcal N(h,1)$ with $h\ge0$: all information about
local depletion is summarized by one normally distributed statistic whose
mean is the depletion strength.
\end{theorem}

The theorem makes precise the effect of a common nuisance level.  The mean of
$j^q$ is removed, and only the centered contrast $g_{j,q,m}$ contributes
information.  ``Matched'' in
\eqref{eq:matchedalternative} means that the common nuisance parameter has
been adjusted so that the average parameter remains $\eta_0$.  Because this
common level is unknown under both hypotheses, only the centered departures
can distinguish depletion from the limiting model.

\begin{remark}[Terminology]
The quadratic expansion \eqref{eq:profilelan} is commonly called profiled local
asymptotic normality, and the family
$\{\mathcal N(h,1):h\ge0\}$ is a one-sided Gaussian shift experiment
\citep{lecam2000,vdvaart1998}.  These terms provide useful context, but the
proof below uses only Taylor's theorem, the law of large numbers, and the
Lindeberg central limit theorem.
\end{remark}

The one-sidedness in the limit comes from the physical constraint $\rho\ge0$,
not from profiling.  If the first coefficient is $+c$ rather than $-c$, one
reverses the sign of $\Delta_m$; the Gaussian shift and the half-$\chi_1^2$
null calibration below are unchanged.

\begin{corollary}[Critical observation scale]
\label{cor:critical}
Set $\rho=N^{-1}$.  Then
\begin{equation}
 h_m
 =
 c\{I(\eta_0)v_q\}^{1/2}
 \frac{m^{q+1/2}}{N^q}\{1+o(1)\}.
\label{eq:hasym}
\end{equation}
Thus the nondegenerate finite-versus-infinite boundary is
\begin{equation}
 {\quad m\asymp N^{2q/(2q+1)}.\quad}
\label{eq:boundary}
\end{equation}
\end{corollary}

\begin{corollary}[Exact first-order cost of an unknown level]
\label{cor:profilingloss}
If the common level were known, the leading information would be
\[
 c^2I(\eta_0)\rho^{2q}\sum_{j<m}j^{2q}
 =\frac{c^2I(\eta_0)}{2q+1}\rho^{2q}m^{2q+1}\{1+o(1)\}.
\]
After profiling, the leading constant is $c^2I(\eta_0)v_q$.  Thus profiling
retains exactly the fraction
\begin{equation}
 \frac{q^2}{(q+1)^2}
\label{eq:profilingloss}
\end{equation}
of the known-level information, or the fraction $q/(q+1)$ of the local
parameter $h$.  It changes the sharp constant but not the exponent.  For JM,
the local parameter is reduced by exactly one half.
\end{corollary}

\begin{corollary}[Information trichotomy]
\label{cor:trichotomy}
Suppose $m\rho\to0$.
\begin{enumerate}
\item If $h_m\to0$, there is a homogeneous nuisance value $\eta_m^\star$ such
that
\[
 \KL\!\left(
   \bigotimes_{j<m}P_{\eta+a(\rho j)},
   P_{\eta_m^\star}^{\otimes m}
 \right)\longrightarrow0.
\]
In particular, these two product laws merge in total variation.
\item If $h_m\to h\in(0,\infty)$, the the limiting profile-likelihood contrast is that of
$\mathcal N(h,1)$ versus $\mathcal N(0,1)$.
\item If $h_m\to\infty$, the centered-score test separates the finite and
homogeneous models consistently.
\end{enumerate}
\end{corollary}

The central sequence in \eqref{eq:central} contains the unknown $\eta_0$.
Let $\widehat\eta_0$ maximize the homogeneous likelihood and define the
feasible residual-score statistic
\begin{equation}
 \widehat\Delta_m
 =-\frac{\sum_{j<m}g_{j,q,m}
 \dot\ell_{\widehat\eta_0}(Y_{j,m})}
 {\{I(\widehat\eta_0)V_{q,m}\}^{1/2}}.
\label{eq:feasiblescore}
\end{equation}
Root-$m$ consistency, a score expansion, and
$\sum_jg_{j,q,m}=0$ give
\begin{equation}
 \widehat\Delta_m=\Delta_m+o_{\Pp}(1)
\label{eq:plugin}
\end{equation}
under the limiting model and the matched local alternatives; details are included
in Appendix~\ref{app:lan}.

\begin{proposition}[Local power envelope]
\label{prop:power}
At the critical scale, the test that rejects homogeneity for
$\widehat\Delta_m>z_{1-\alpha}$ has asymptotic size $\alpha$ and limiting power
\begin{equation}
  \beta_\alpha(h)=1-\Phi(z_{1-\alpha}-h).
\label{eq:power}
\end{equation}
No other sequence of tests with asymptotic size at most $\alpha$ uniformly in
the local nuisance value has larger pointwise limiting power against the
matched alternative \eqref{eq:matchedalternative}.
\end{proposition}

\begin{proof}
First reveal \(\eta_0\).  For each \(m\), the Neyman-Pearson lemma then
bounds the power of every test in the original composite problem, because knowing the nuisance value cannot lower the oracle envelope.  Let \(L_m\) denote the likelihood ratio of the matched alternative to the resulting
limiting-model law.  Under this law, the Taylor calculation in the proof of
Theorem~\ref{thm:lan} gives
\[
 \log L_m
 =
 h\Delta_m-\frac{h^2}{2}+o_{\Pp}(1)
 \ \Longrightarrow\
 \mathcal N(-h^2/2,h^2).
\]
Repeating the same Taylor expansion under the matched alternative (with the
uniform remainder bounds unchanged) and using
\eqref{eq:alternativeclt} gives
$\log L_m\Rightarrow\mathcal N(h^2/2,h^2)$ directly.
The limiting null distribution is continuous, so the size-$\alpha$ likelihood
ratio cutoff converges to the cutoff $\Delta=z_{1-\alpha}$.  Its limiting
power is therefore \eqref{eq:power}.  By \eqref{eq:plugin}, the feasible test
attains the oracle bound and hence the composite envelope.
\end{proof}


\section{Contact relative to a nuisance tangent space}
\label{sec:tangent}

The preceding theorem centers a polynomial because the only nuisance
direction is the constant vector.  With several nuisance directions, however,
``orthogonal'' means orthogonal in the local Fisher metric, not necessarily in
ordinary Euclidean distance.  This distinction matters whenever the limiting-model
parameter, and hence its Fisher information, varies over the design. Here we present the more general and subsuming theorem (Theorem \ref{thm:tangent}).

Let $p\ge1$ be a fixed integer, the dimension of the nuisance trend, and
let
\[
 z(t)=\{1,z_1(t),\ldots,z_{p-1}(t)\}^\mathsf T
\]
collect $p$ continuous, linearly independent functions on $[0,1]$, the
first identically one.  Write $t_{j,m}=j/m$ and consider
\begin{equation}
 Y_{j,m}\sim P_{z(t_{j,m})^\mathsf T\beta+a(\rho j)},
 \qquad j=0,\ldots,m-1,
\label{eq:tangentmodel}
\end{equation}
where $\beta\in\R^p$ is an unknown coefficient vector.  The
limiting model is obtained by setting  $\rho=0$ with $\beta$ remaining free; Section~\ref{sec:general}
is recovered at $p=1$, $\beta=\eta$.  The first coordinate of $\beta$
plays the role of the common level, and the remaining coordinates model
known-shape, unknown-magnitude trends such as the effort drift of
Section~\ref{sec:seconddrift} below.
Write $z_{j,m}=z(t_{j,m})$ and let $Z_m$ be the matrix with rows
$z_{j,m}^\mathsf T$.  The baseline parameter under the limiting model is allowed to vary with $m$:
\begin{equation}
 \theta_{j,0}=z_{j,m}^\mathsf T\beta_{0,m},
 \qquad
 D_m=\operatorname{diag}\{I(\theta_{0,0}),\ldots,I(\theta_{m-1,0})\}.
\label{eq:fisherweights}
\end{equation}

We call \(\mathbb R^m\), equipped with
\(\langle v,w\rangle_{m,I}=v^\mathsf TD_mw\), the
\emph{Fisher-weighted design-shift space}.  A vector \(v\) records the
local parameter shifts across the \(m\) design positions, while
\(\operatorname{col}(Z_m)\) is the nuisance design-shift subspace.

For $x_{k,m}=(0^k,1^k,\ldots,(m-1)^k)^\mathsf T$, define the
$D_m$-orthogonal projection and residual by
\begin{equation}
 \Pi_{m,I}=Z_m(Z_m^\mathsf TD_mZ_m)^{-1}Z_m^\mathsf TD_m,
 \qquad
 R_{k,m}^{I}=(I_m-\Pi_{m,I})x_{k,m},
 \qquad
 \mathcal J_{k,m}=(R_{k,m}^{I})^\mathsf TD_mR_{k,m}^{I}.
\label{eq:tangentresidual}
\end{equation}
Thus $Z_m^\mathsf TD_mR_{k,m}^{I}=0$.  The superscript $I$ emphasizes
that the residual is determined by the information geometry.  When
$I(\theta_{j,0})$ is constant, this reduces to ordinary least-squares
projection.

\begin{definition}[Effective contact order]
\label{def:effective}
Let $a(x)=\sum_{k\ge1}a_kx^k$ near zero.  The effective contact order is
the first $q$ for which $a_q\ne0$ and $R_{q,m}^{I}\ne0$, provided that every
active lower-order trend is in $\operatorname{col}(Z_m)$.  We assume, after
orienting the alternative, that $a_q=-c<0$, that
$\mathcal J_{q,m}/m^{2q+1}$ converges to a positive limit, and that
\begin{equation}
 \max_j
 \frac{I(\theta_{j,0})^{1/2}|R_{q,m,j}^{I}|}
      {\sqrt{\mathcal J_{q,m}}}
 =O(m^{-1/2}).
\label{eq:effectiveconditions}
\end{equation}
\end{definition}

This definition is the observable version of analytic contact: a
Taylor coefficient counts only when it points outside the nuisance tangent
space.  Membership of the lower-order trends in
$\operatorname{col}(Z_m)$ makes their absorption exact under any positive
weighting: a vector of $\operatorname{col}(Z_m)$ is annihilated by
$I_m-\Pi_{m,I}$ regardless of $D_m$, so Fisher weighting changes the
residual and its norm, never which terms are removed.

The maximal-weight condition in Definition~\ref{def:effective} is a
Lindeberg condition: no single observation may carry a nonvanishing fraction
of the projected information.  The following lemma shows that it is automatic
under the usual smooth fixed-design assumptions.  It also identifies the
limit of $\mathcal J_{q,m}/m^{2q+1}$ directly from the limiting weighted
regression.

\begin{lemma}[Smooth-design verification]
\label{lem:smoothdesign}
Suppose $z$ is continuous and bounded on $[0,1]$, $\beta_{0,m}\to\beta_0$,
and
\[
 w_m(t)=I\{z(t)^\mathsf T\beta_{0,m}\}
 \longrightarrow
 w(t)=I\{z(t)^\mathsf T\beta_0\}
\]
uniformly, where $w$ is positive and continuous.  Assume
\[
 G_I=\int_0^1w(t)z(t)z(t)^\mathsf Tdt>0,
 \qquad
 b_q=G_I^{-1}\int_0^1w(t)z(t)t^qdt,
\]
and put $\psi_q(t)=t^q-z(t)^\mathsf Tb_q$.  If $\psi_q$ is not
identically zero, then
\begin{align}
 \max_{j<m}
 \left|m^{-q}R_{q,m,j}^{I}-\psi_q(t_{j,m})\right|&\to0,
 \label{eq:residualcontinuum}\\
 \mathcal J_{q,m}
 &=m^{2q+1}\left\{\int_0^1w(t)\psi_q(t)^2dt+o(1)\right\}.
 \label{eq:Jcontinuum}
\end{align}
Consequently the conditions in \eqref{eq:effectiveconditions} hold, and
\[
 \lim_{m\to\infty}\frac{\mathcal J_{q,m}}{m^{2q+1}}
 =\int_0^1w(t)\psi_q(t)^2dt>0.
\]
\end{lemma}

\begin{proof}
Let
\[
 A_m=m^{-1}Z_m^\mathsf TD_mZ_m,
 \qquad
 c_{q,m}=m^{-1}\sum_{j<m}w_m(t_{j,m})z(t_{j,m})t_{j,m}^q.
\]
Uniform convergence and Riemann sums give
$A_m\to G_I$ and
$c_{q,m}\to\int_0^1w(t)z(t)t^qdt$.  The coefficient in the
weighted projection of $x_{q,m}$, divided by $m^q$, is
$A_m^{-1}c_{q,m}\to b_q$.  Boundedness of $z$ now gives
\eqref{eq:residualcontinuum} uniformly in $j$.  A second Riemann-sum
argument yields \eqref{eq:Jcontinuum}.  Finally, positivity and
boundedness of $w_m$, uniform boundedness of $m^{-q}R_{q,m,j}^{I}$, and
\eqref{eq:Jcontinuum} imply
\[
 \max_{j<m}
 \frac{I(\theta_{j,0})^{1/2}|R_{q,m,j}^{I}|}
      {\sqrt{\mathcal J_{q,m}}}
 =O\left(\frac{m^q}{m^{q+1/2}}\right)
 =O(m^{-1/2}).
\]
\end{proof}

Nothing in the proof uses integrality of $q$: for any real $r>0$ the map
$t\mapsto t^r$ is continuous on $[0,1]$, so the lemma verifies the two
conditions required by Corollary~\ref{cor:powercontact} below verbatim,
with $\psi_r(t)=t^r-z(t)^\mathsf Tb_r$.

Put $a_m(\rho)=\{a(\rho j):j<m\}$ and
\begin{equation}
 d_m(\rho)=(I_m-\Pi_{m,I})a_m(\rho),
 \qquad
 \gamma_m(\rho)
 =(Z_m^\mathsf TD_mZ_m)^{-1}Z_m^\mathsf TD_ma_m(\rho).
\label{eq:projecteddepletion}
\end{equation}
Since $\Pi_{m,I}a_m(\rho)=Z_m\gamma_m(\rho)$, a finite model with nuisance
$\beta$ has, to first order in the local Fisher metric, closest limiting-model
parameter
\begin{equation}
 \beta_m^\star=\beta+\gamma_m(\rho).
\label{eq:closestbeta}
\end{equation}
Conversely, the alternative matched to the limiting-model baseline $\beta_{0,m}$
uses nuisance value $\beta_{0,m}-\gamma_m(\rho)$ and has parameter vector
$Z_m\beta_{0,m}+d_m(\rho)$.  This distinction permits nuisance shifts such as
$m/N$: they vanish, but need not be of root-$m$ order.  The
closest-limiting-model identity \eqref{eq:closestbeta} is also the matched value
at which the merging statement is computed: the Kullback-Leibler expansion
closing Appendix~\ref{app:tangent} compares each finite law with the
limiting-model at exactly $\beta_m^\star$, so the trichotomy is a statement
about Fisher-closest pairs.

For this section write
\[
 \pl_m(\rho)=
 \sup_\beta\sum_{j<m}
 \ell_{z_{j,m}^\mathsf T\beta+a(\rho j)}(Y_{j,m}).
\]

\begin{theorem}[Profile likelihood after nuisance projection]
\label{thm:tangent}
Suppose (S1)-(S2) hold uniformly on a fixed neighborhood of all
$\theta_{j,0}$ and the local shifts below.  Let the regressors be uniformly
bounded, let $\beta_{0,m}\to\beta_0$, and suppose
\begin{equation}
 m^{-1}Z_m^\mathsf TD_mZ_m\longrightarrow G_I>0.
\label{eq:weightedgram}
\end{equation}
Suppose also that Definition~\ref{def:effective} holds, that the likelihood under the limiting model is concave, and that its MLE satisfies
$\widehat\beta_0-\beta_{0,m}=O_{\Pp}(m^{-1/2})$.  Define
\begin{equation}
 h_m=c\rho^q\sqrt{\mathcal J_{q,m}},
 \qquad
 \Delta_m^{\perp}
 =-\frac{\sum_{j<m}R_{q,m,j}^{I}
 \dot\ell_{\theta_{j,0}}(Y_{j,m})}
 {\sqrt{\mathcal J_{q,m}}}.
\label{eq:tangentcentral}
\end{equation}
Then, with
$\rho_m(h)=\{h/(c\sqrt{\mathcal J_{q,m}})\}^{1/q}$,
under the limiting-model law $\bigotimes_{j<m}P_{\theta_{j,0}}$,
\begin{equation}
 \pl_m\{\rho_m(h)\}-\pl_m(0)
 =h\Delta_m^\perp-\frac{h^2}{2}+o_{\Pp}(1)
\label{eq:tangentlan}
\end{equation}
uniformly for $h$ in compact subsets of $[0,\infty)$, and
$\Delta_m^\perp\Rightarrow\mathcal N(0,1)$.  Under the matched alternative
$Y_{j,m}\sim P_{\theta_{j,0}+d_{m,j}\{\rho_m(h_0)\}}$,
$\Delta_m^\perp\Rightarrow\mathcal N(h_0,1)$.

The result extends to a feasible statistic.  Set
\[
 \widehat\theta_{j,0}=z(t_{j,m})^\mathsf T\widehat\beta_0,
 \qquad
 \widehat D_m=\operatorname{diag}\{I(\widehat\theta_{j,0}):j<m\},
\]
\[
 \widehat\Pi_{m,I}
 =Z_m(Z_m^\mathsf T\widehat D_mZ_m)^{-1}Z_m^\mathsf T\widehat D_m,
 \quad
 \widehat R_{q,m}^{I}=(I_m-\widehat\Pi_{m,I})x_{q,m},
 \quad
 \widehat{\mathcal J}_{q,m}
 =(\widehat R_{q,m}^{I})^\mathsf T\widehat D_m
   \widehat R_{q,m}^{I}.
\]
Then
\begin{equation}
 \widehat\Delta_m^\perp
 =-\frac{\sum_{j<m}\widehat R_{q,m,j}^{I}
 \dot\ell_{\widehat\theta_{j,0}}(Y_{j,m})}
 {\sqrt{\widehat{\mathcal J}_{q,m}}}
 =\Delta_m^\perp+o_{\Pp}(1).
\label{eq:tangentfeasible}
\end{equation}
Consequently, the merge-critical-separate trichotomy and the Gaussian power
formula in \eqref{eq:power} hold with the effective local parameter $h_m$
defined above.  If $\rho=N^{-1}$, the critical scale remains
$m\asymp N^{2q/(2q+1)}$.
\end{theorem}

\begin{corollary}[Local power contact beyond analyticity]
\label{cor:powercontact}
Analyticity of $a$ in Theorem~\ref{thm:tangent} may be replaced by the
one-sided local expansion
\begin{equation}
 a(x)=\sum_{\ell=1}^{\nu}b_\ell x^{r_\ell}-cx^r+o(x^r),
 \qquad
 0<r_1<\cdots<r_\nu<r,
 \qquad c>0,
 \label{eq:powercontact}
\end{equation}
where $r>0$ need not be an integer.  For
$x_{u,m}=(0^u,1^u,\ldots,(m-1)^u)^\mathsf T$, suppose every
$x_{r_\ell,m}$ belongs to $\operatorname{col}(Z_m)$.  Define
\[
 R_{r,m}^{I}=(I_m-\Pi_{m,I})x_{r,m},
 \qquad
 \mathcal J_{r,m}=(R_{r,m}^{I})^\mathsf TD_mR_{r,m}^{I},
\]
and suppose the two conditions in \eqref{eq:effectiveconditions} hold
with $q$ replaced by $r$.  Then all conclusions of
Theorem~\ref{thm:tangent} hold with
\[
 h_m=c\rho^r\sqrt{\mathcal J_{r,m}},
 \qquad
 \rho_m(h)=\left\{\frac{h}{c\sqrt{\mathcal J_{r,m}}}\right\}^{1/r}.
\]
In particular, if $\rho=N^{-1}$ and
$\mathcal J_{r,m}\asymp m^{2r+1}$, the critical scale is
\[
 m\asymp N^{2r/(2r+1)}.
\]
\end{corollary}

\begin{proof}
Write
\[
 \varepsilon(x):=a(x)-\sum_{\ell=1}^{\nu}b_\ell x^{r_\ell}+cx^r,
\]
the remainder of the expansion \eqref{eq:powercontact}, so that
$\varepsilon(x)=o(x^r)$ as $x\downarrow0$.  Projection removes every term
with exponent $r_\ell$.  Because $m\rho\to0$ on the displayed local
scale,
\[
 \max_{j<m}|\varepsilon(\rho j)|
 =o\{(m\rho)^r\}.
\]
Orthogonal projection is a contraction in the $D_m$ norm, while the
information weights are bounded above and away from zero.  Hence the
projected remainder has $D_m$ norm
\[
 o\{\rho^rm^{r+1/2}\}
 =o\{\rho^r\sqrt{\mathcal J_{r,m}}\}.
\]
The leading projected displacement is therefore
$-c\rho^rR_{r,m}^{I}$, uniformly on compact $h$-sets.  The likelihood
Taylor expansion, maximal-weight central limit theorem, and nuisance
orthogonality argument in the proof of Theorem~\ref{thm:tangent} now
apply verbatim with $q$ replaced by $r$.
\end{proof}

Thus a cusp such as $a(x)=-cx^r+o(x^r)$ has contact order $r$, even when
$r$ is noninteger\footnote{For example,  consider  depletion with  
 $N$ units having (depletion) rates $\theta W_1,\ldots,\theta W_N$, where the common
scale $\theta>0$ is unknown and
$
 W_1,\ldots,W_N\stackrel{\mathrm{iid}}{\sim}F
$ 
. When $F$ is heavy-tailed, we conjecture that such fractional boundary governs this  removal
experiment;  proving it, however is beyond the scope of this paper.}.  Analyticity supplies a convenient way to locate the first
surviving term; the likelihood boundary itself depends only on the local power
and its residual Fisher norm.

The proof, in Appendix~\ref{app:tangent}, also explains why the projection
must be taken before assigning a contact order.  The deterministic
orthogonality $Z_m^\mathsf TD_mR_{q,m}^{I}=0$ eliminates the leading nuisance
cross information; the corresponding observed cross term is $o_{\Pp}(1)$,
not identically zero.  Fisher weighting generally changes the sharp constant.
By contrast, enlarging the nuisance space can remove a complete Taylor term,
change $q$, and thereby change the boundary exponent.

\begin{remark}
\label{rem:two-scales}
At the critical scale, two notions of size behave differently as $m$ grows.
The parameter shifts that remain after nuisance projection have
root-mean-square magnitude of order $m^{-1/2}$, irrespective of the contact
order $q$. By contrast, the design samples the curve over an argument interval whose width is of scale $m^{-1/(2q)}$.

\emph{(i) Nuisance-orthogonal shifts.}
At $\rho=\rho_m(h)$, for fixed $h>0$, the residual departure
$d_m=d_m\{\rho_m(h)\}$ satisfies
\[
    d_m^{\mathsf T}D_m d_m=h^2\{1+o(1)\}
\]
(Appendix~\ref{app:tangent}). Since $I$ is continuous, positive, and bounded
along the baseline curve,
\[
    \sum_{j<m} d_{m,j}^2 \asymp h^2,
\]
and hence the root-mean-square shift satisfies
\[
  \left(m^{-1}\sum_{j<m}d_{m,j}^2\right)^{1/2}
  \asymp h\,m^{-1/2}.
\]
Moreover, under the maximal-weight condition \eqref{eq:effectiveconditions},
\[
    \max_{j<m}|d_{m,j}|=O(m^{-1/2}).
\]
Neither scale depends on $q$: at a nondegenerate boundary the residual
Fisher information is of constant order and, under the maximal-weight
condition, is spread across the expanding design at root-$m$ scale. This root-$m$ residual scale is the analogue, for the present non-identically distributed experiment, of the usual Pitman localization.

\emph{(ii) Sampled range of the curve.}
The design evaluates $a$ only at the points $\rho j$, $j<m$, and hence over
the interval $[0,\rho(m-1)]$. Since
\[
    h_m=c\rho^q\sqrt{J_{q,m}},
    \qquad
    J_{q,m}\asymp m^{2q+1},
\]
the relation $h_m\asymp1$ is equivalent to
\[
    (\rho m)^q\asymp m^{-1/2},
\]
or
\begin{equation}\label{eq:sampled-range}
    \rho m\asymp m^{-1/(2q)}.
\end{equation}
Thus, at the detection boundary, the observable parameter shifts remaining after nuisance projection are always of root-$m$ size, irrespective of $q$ (this is the analogue here of the usual Pitman localization), while the range of the curve argument that must be sampled to
produce it depends on $q$. Higher contact order makes this range shrink more
slowly: a flatter observable departure requires the design to extend farther
along the curve before a root-$m$ residual signal emerges.

With $\rho=N^{-1}$, \eqref{eq:sampled-range} is equivalently
\[
    m\asymp N^{2q/(2q+1)}.
\]
Indeed, at this scale both $m/N$ and $m^{-1/(2q)}$ are of order
$N^{-1/(2q+1)}$, and $m/N$ is the fraction of the population depleted.
\end{remark}

\begin{corollary}[Polynomial nuisance and a closed information constant]
\label{cor:legendre}
Suppose $I(\theta_{j,0})\equiv I_0$ and
$\operatorname{col}(Z_m)=\operatorname{span}\{1,j,\ldots,j^{q-1}\}$.
Write $\mathcal P_{q-1}$ for the polynomials on $[0,1]$ of degree at most
$q-1$ and $\operatorname{proj}_{\mathcal P_{q-1}}$ for orthogonal
projection onto it in $L^2(0,1)$, so that
$\operatorname{proj}_{\mathcal P_{q-1}}t^q$ is the best mean-square
approximation of $t^q$ by a polynomial of degree below $q$.
Then
\begin{equation}
 \mathcal J_{q,m}
 =I_0w_qm^{2q+1}\{1+o(1)\},
 \qquad
 w_q
 =\int_0^1
 \left\{t^q-\operatorname{proj}_{\mathcal P_{q-1}}t^q\right\}^2dt
 =\frac{(q!)^4}{[(2q)!]^2(2q+1)}.
\label{eq:legendreconstant}
\end{equation}
Relative to known nuisance coefficients, profiling retains the exact
first-order information fraction
\begin{equation}
 (2q+1)w_q=\binom{2q}{q}^{-2}.
\label{eq:polynomialloss}
\end{equation}
\end{corollary}

\begin{proof}
Let $L_q(t)=P_q(2t-1)$ be the shifted Legendre polynomial.  It is orthogonal
to $\mathcal P_{q-1}$ on $[0,1]$, has leading coefficient
$\binom{2q}{q}$, and satisfies
$\int_0^1L_q(t)^2dt=(2q+1)^{-1}$.  Hence the residual of $t^q$ is
$L_q(t)/\binom{2q}{q}$.  A Riemann-sum approximation proves
\eqref{eq:legendreconstant}; division by
$\sum_{j<m}j^{2q}\sim m^{2q+1}/(2q+1)$ proves
\eqref{eq:polynomialloss}.
\end{proof}

\begin{remark}
\label{rem:polynomialhypothesis}
The assumption
$\operatorname{col}(Z_m)=\operatorname{span}\{1,j,\ldots,j^{q-1}\}$ is
stronger than Definition~\ref{def:effective} 
and should not be read as a restatement of it.  The definition requires
only that the \emph{active} Taylor coefficients below order $q$ lie in
$\operatorname{col}(Z_m)$; it permits a nuisance space containing fewer
monomials, or none at all.  The constant $w_q$ is the value taken when the
nuisance space is exactly $\mathcal P_{q-1}$.

Two departures illustrate the difference.  If $a_1=\cdots=a_{q-1}=0$, so
that the profile has genuine $q$th-order tangency and no lower-order trend
requires absorption, then an unknown common level alone already gives
effective contact order $q$.  The relevant residual is then $j^q$ centered
rather than residualized on $\operatorname{span}\{1,\ldots,j^{q-1}\}$, and
the constant is $v_q$ of \eqref{eq:vdef} rather than $w_q$; for $q=2$ these
differ by a factor of sixteen, since $v_2=4/45$ and $w_2=1/180$.
Conversely, a nuisance space containing non-polynomial regressors, or
monomials of degree exceeding $q$, yields neither constant, and
Lemma~\ref{lem:smoothdesign} 
must be used to evaluate $\lim_m\mathcal J_{q,m}/m^{2q+1}$ directly.

Corollaries~\ref{cor:profilingloss} and~\ref{cor:legendre} therefore
bracket the polynomial range: the former treats the smallest nuisance
space compatible with contact order $q$, retaining the fraction
$q^2/(q+1)^2$, and the latter the largest, retaining
$\binom{2q}{q}^{-2}$.  The two agree at $q=1$, where both nuisance spaces
reduce to the constants and both fractions equal $1/4$.  In either case
the exponent $2q/(2q+1)$ is unchanged: the nuisance space determines
\emph{which} Taylor coefficient is the first observable one, and then
charges a constant for profiling it.  Explicitly, for $q=2$ one has
$\operatorname{proj}_{\mathcal P_1}t^2=t-1/6$, so the residual is
$t^2-t+1/6$, which is the continuum limit of
\eqref{eq:quadraticresidual} after division by $m^2$.
\end{remark}

\begin{remark} Magnitude and decay of the information loss: 
\label{rem:brutal}
It is interesting to examine the fraction of retained information,
$(2q+1)w_q=\binom{2q}{q}^{-2}$, whose decay in $q$ is severe: for small $q$  the table below shows a fast  decay 
\begin{center}
\begin{tabular}{@{}lccccc@{}}
\toprule
$q$ & 1 & 2 & 3 & 4 & 5\\
\midrule
retained information $\binom{2q}{q}^{-2}$
 & $1/4$ & $1/36$ & $1/400$ & $1/4900$ & $1/63504$\\
retained local parameter $\binom{2q}{q}^{-1}$
 & $1/2$ & $1/6$ & $1/20$ & $1/70$ & $1/252$\\
sample-size penalty $\binom{2q}{q}^{2/(2q+1)}$
 & $1.59$ & $2.05$ & $2.35$ & $2.57$ & $2.73$\\
\bottomrule
\end{tabular}
\end{center}
and for large values of $q$, since $\binom{2q}{q}\sim4^q/\sqrt{\pi q}$, the retained fraction decays as
$\pi q\,16^{-q}$: each additional polynomial nuisance direction costs
roughly a further factor of sixteen in information. 

The penalty is nevertheless confined to the constant (while the exponent $2q/(2q+1)$ is untouched).  Thus, because
$h_m\asymp m^{q+1/2}N^{-q}$, an information loss of factor
$\binom{2q}{q}^{-2}$ inflates the number of observations required for a
given local parameter by $\binom{2q}{q}^{2/(2q+1)}$. Moreover, this sample-size inflation grows slowly,
  under three
for every $q\le5$ and approaches monotonically to $4$.  Profiling a
richer polynomial nuisance space is therefore expensive but never fatal,
whereas enlarging that space is what changes $q$ itself.
\end{remark}


\subsection{A second-order example}
\label{sec:seconddrift}

Suppose that failure-detection effort changes during a study.  A simple
exponential model allowing unknown log-linear effort drift is
\begin{equation}
 X_{j,N}\sim\Exp\!\left[
 e^{\eta+\gamma(j/m)}\left(1-\frac jN\right)\right],
 \qquad j=0,\ldots,m-1,
\label{eq:driftmodel}
\end{equation}
where the limiting-model omits $1-j/N$ but retains both $\eta$ and $\gamma$.
For the log-rate parameter $\theta=\log\lambda$, the score is
$1-\lambda X$ and $I(\theta)=1$.  Thus the Fisher projection in this example
is exactly the Euclidean projection.
On the log-rate scale,
\[
 \log(1-j/N)=-\frac jN-\frac{j^2}{2N^2}
 +O\{(j/N)^3\}.
\]
The first term is indistinguishable from changing $\gamma$ by $-m/N$.
The first surviving trend is quadratic.  Its exact residual on
$\operatorname{span}\{1,j\}$ is
\begin{equation}
 R_{2,m,j}^{I}=j^2-(m-1)j+\frac{(m-1)(m-2)}6,
 \qquad
 \mathcal J_{2,m}=\frac{m(m^2-1)(m^2-4)}{180}.
\label{eq:quadraticresidual}
\end{equation}
Thus \eqref{eq:driftmodel} has effective contact $q=2$, $c=1/2$, and
\begin{equation}
 h_{m,N}=\frac{1}{2N^2}\sqrt{\mathcal J_{2,m}}
 \sim\frac{m^{5/2}}{\sqrt{720}\,N^2},
 \qquad
 {m\asymp N^{4/5}}.
 \label{eq:fourfifths}
\end{equation}
The statistic \eqref{eq:tangentfeasible} uses the limiting-model MLEs of
$\eta$ and $\gamma$ and log-rate scores
$1-\widehat\lambda_jX_{j,N}$.  This example demonstrate   two separate nuisance effects:
an unknown level changes only a constant, whereas an unknown trend can absorb
a whole Taylor coefficient and change the detection exponent.

\section{Exponential rate arrays: an exact check}
\label{sec:exp}
Exponential rate arrays are one family in which the profile likelihood
admits a closed form.  This yields the likelihood-ratio law directly,
independently of Theorem~\ref{thm:lan}, and supplies the exact
finite-sample quantities used for the Jelinski-Moranda consequences of
Section~\ref{sec:jm} and the simulations of Section~\ref{sec:simulation}.

Consider independent variables
\begin{equation}
 X_{j,N}\sim
 \Exp\{\alpha f(j/N)\},\qquad j=0,\ldots,m-1,
\label{eq:exparray}
\end{equation}
where $\alpha>0$ is unknown, $f(0)=1$, and
\begin{equation}
  f(x)=1-cx^q+O(x^{q+1}).
\label{eq:fcontact}
\end{equation}
The limiting model consists of iid $\Exp(\beta)$ variables with unknown
$\beta>0$.  Taking $\theta=\log\alpha$ and $a(x)=\log f(x)$ embeds
\eqref{eq:exparray} in \eqref{eq:analyticmodel}.  The log-rate score is
$1-\alpha X$ and its Fisher information is one.

Let $\Lambda_{m,N}$ be the maximized log likelihood under
\eqref{eq:exparray} minus that under the limiting model.  There is also a
useful exact representation.  Put
\[
 b_{j,N}=f(j/N)^{-1},\qquad
 \bar b_{m,N}=\frac1m\sum_{j=0}^{m-1}b_{j,N},\qquad
 \sigma_{j,N}=\frac{b_{j,N}}{\bar b_{m,N}}.
\]
Under the finite model,
\[
 E_j=\alpha f(j/N)X_{j,N}\stackrel{\mathrm{iid}}{\sim}\Exp(1).
\]

\begin{proposition}[Exact profile identity]
\label{prop:exactexp}
Under \eqref{eq:exparray},
\begin{equation}
 \Lambda_{m,N}
 =
 -\sum_{j=0}^{m-1}\log \sigma_{j,N}
 +
 m\log
 \left(
  \frac{\sum_{j=0}^{m-1}\sigma_{j,N}E_j}
      {\sum_{j=0}^{m-1}E_j}
 \right).
\label{eq:exactprofile}
\end{equation}
\end{proposition}

The identity, however, relies on the nuisance being a pure scale: profiling $\alpha$
out of an exponential likelihood is explicit, and the factor
$\bar b_{m,N}$ cancels.  No such closed form is available for the
two-dimensional nuisance of Section~\ref{sec:seconddrift}, so the general
theorems are required rather than convenient.

Define
\begin{equation}
 \kappa_{m,N}
 =
 c^2v_q\frac{m^{2q+1}}{N^{2q}}.
\label{eq:kappa}
\end{equation}

\begin{corollary}[Likelihood-ratio limit]
\label{cor:lrlimit}
If $m/N\to0$ and $\kappa_{m,N}\to\kappa\in[0,\infty)$, then
\begin{align}
 \Lambda_{m,N}
 &\Longrightarrow
 \mathcal N(\kappa/2,\kappa)
 &&\text{under the finite model},\label{eq:lrfinite}\\
 \Lambda_{m,N}
 &\Longrightarrow
 \mathcal N(-\kappa/2,\kappa)
 &&\text{under the limiting-model}.
\label{eq:lrnull}
\end{align}
\end{corollary}

The means $\pm\kappa/2$ and common variance $\kappa$ are the canonical
log-likelihood-ratio form of a Gaussian shift experiment.  In particular,
$\kappa=h^2$.

\section{Sharp JM consequences at the sparse boundary}
\label{sec:jm}

For JM, $f(x)=1-x$, so $q=c=1$.  The exact quantities
$V_m$, $\Delta_m$, and $h_{m,N}$ were introduced in
\eqref{eq:Vm}-\eqref{eq:hjm}.  The general theorem now makes the motivating
calculation \eqref{eq:jmheuristic} rigorous and supplies its inferential
consequences.

\begin{theorem}[JM boundary experiment]
\label{thm:jm}
Let $m\to\infty$, $N\to\infty$, and $m/N\to0$.
\begin{enumerate}
\item If $h_{m,N}\to h\in[0,\infty)$, then under JM,
\[
  \Delta_m\Longrightarrow\mathcal N(h,1),
\]
and the profiled likelihood process converges on compact subsets of
$[0,\infty)$ to
\[
  u\longmapsto u\Delta-\frac{u^2}{2},
  \qquad \Delta\sim\mathcal N(h,1).
\]
\item If $h_{m,N}\to0$, the JM experiment and an iid exponential limiting-model
with a suitably matched rate merge in total variation.
\item If $h_{m,N}\to\infty$, depletion can be detected consistently.
\end{enumerate}
\end{theorem}

Since $h_{m,N}$ is of constant order exactly when $m\asymp N^{2/3}$,
Theorem~\ref{thm:jm} identifies the first scale at which a finite 
population is statistically visible.

The following corollary is partially known already (its eq. (\ref{eq:halfexact})), but we bring it here for completeness.

\begin{corollary}[Probability of an infinite estimate]
\label{cor:pinfty}
If $h_{m,N}\to h$, then
\begin{equation}
  \Pp_N(\widehat N=\infty)
  \longrightarrow
  \Phi(-h).
\label{eq:pinfty}
\end{equation}
Under the exact iid limiting model,
\begin{equation}
  \Pp_\infty(\widehat N=\infty)=\frac12
\label{eq:halfexact}
\end{equation}
for every $m\ge2$.
\end{corollary}

\begin{proof}
The first assertion follows from Proposition~\ref{prop:existence} and
Theorem~\ref{thm:jm}.  Under the limiting model, reversal of the exchangeable
vector $(X_0,\ldots,X_{m-1})$ maps $W_m$ to $m-1-W_m$.  Since the
distribution is continuous, it assigns equal probabilities to the two sides
of $\bar j_m=(m-1)/2$.
\end{proof}

For the explicit critical sequence $m=sN^{2/3}\{1+o(1)\}$,
\begin{equation}
 h=\frac{s^{3/2}}{\sqrt{12}},
 \qquad
 \Pp_N(\widehat N=\infty)
 \longrightarrow
 \Phi\left(-\frac{s^{3/2}}{\sqrt{12}}\right).
\label{eq:sexplicit}
\end{equation}

\begin{corollary}[Boundary law of the MLE]
\label{cor:mlelaw}
Let
\[
  \widehat h_m=\frac{\sqrt{V_m}}{\widehat N},
  \qquad \Delta_+ := \max(0, \Delta),
  \qquad \frac1\infty:=0, 
  \qquad \frac1 0:=\infty.
\]
If $h_{m,N}\to h$, then
\begin{equation}
  \widehat h_m
  \Longrightarrow
  \Delta_+,
  \qquad
  \Delta\sim\mathcal N(h,1).
\label{eq:hhat}
\end{equation}
Equivalently,
\begin{equation}
 \frac{\widehat N}{\sqrt{V_m}}
 \Longrightarrow
 \frac1{\Delta_+},
\label{eq:Nhat}
\end{equation}
where the limit has an atom of size $\Phi(-h)$ at infinity.
\end{corollary}

\begin{corollary}[Likelihood-ratio law and local power]
\label{cor:lrt}
For testing the limiting-model hypothesis $h=0$ against finite depletion $h>0$,
\begin{equation}
 2\left\{\sup_{N<\infty,\phi}\ell(N,\phi)
 -\sup_{\beta}\ell(\infty,\beta)\right\}
 \Longrightarrow
 (\Delta_+)^2.
\label{eq:lrt}
\end{equation}
Under the limiting-model this has distribution
\begin{equation}
 \frac12\delta_0+\frac12\chi_1^2.
\label{eq:chernoff}
\end{equation}
By Proposition~\ref{prop:power}, the asymptotic power envelope at local
parameter $h$ and level $\alpha$ is
$1-\Phi(z_{1-\alpha}-h)$.
\end{corollary}
\noindent with $\delta_0$ denoting the unit point mass at the origin; i.e.: with
probability one half the constrained maximizer lies at the boundary and the
statistic vanishes identically, and otherwise it is the square of a
positive standard normal variable.

\begin{corollary}[Likelihood confidence set for population size]
\label{cor:confidence}
Let $c_\alpha$ be the $(1-\alpha)$ quantile of $\chi_1^2$ and use the local
profile likelihood to form
\[
 \mathcal C_{h,m}=
 \left\{u\ge0:
 2\left[\sup_{v\ge0}\pl_m(v)-\pl_m(u)\right]\le c_\alpha
 \right\},
\]
where the argument is expressed on the local scale $u=\sqrt{V_m}/N$.
This set has asymptotic coverage at least $1-\alpha$ for every fixed local
$h\ge0$.  In the limiting Gaussian experiment,
\begin{align*}
 u_L(\Delta)&=\max\{0,\Delta-\sqrt{c_\alpha}\},\\
 u_U(\Delta)&=
 \begin{cases}
 \Delta+\sqrt{c_\alpha},&\Delta>0,\\
 \Delta+\sqrt{\Delta^2+c_\alpha},&\Delta\le0.
 \end{cases}
\end{align*}
Then $\mathcal C_h=[u_L(\Delta),u_U(\Delta)]$, whose image for
population size is
\begin{equation}
 \mathcal C_N=
 \begin{cases}
 [\sqrt{V_m}/u_U,\,\infty],&u_L=0,\\
 [\sqrt{V_m}/u_U,\,\sqrt{V_m}/u_L],&u_L>0.
 \end{cases}
\label{eq:Nconfidence}
\end{equation}
\end{corollary}

\begin{proof}
In the limit, $Q(u)=u\Delta-u^2/2$ and its constrained maximizer is
$\widehat u=\Delta_+$.  Solving
$2\{Q(\widehat u)-Q(u)\}\le c_\alpha$ gives the displayed endpoints.
At the true $h$, this likelihood drop is at most $(\Delta-h)^2$; hence its
coverage is at least
$\Pp_h\{(\Delta-h)^2\le c_\alpha\}=1-\alpha$.  The compact-uniform
likelihood expansion transfers the result to the finite experiment.  Since
$\sqrt{c_\alpha}=z_{1-\alpha/2}$, the limiting coverage is exactly
$1-\alpha$ for $h\ge z_{1-\alpha/2}$ and is strictly larger for
$0\le h<z_{1-\alpha/2}$; at the limiting-model  $h=0$ it equals
$1-\alpha/2$.  Indeed, writing $\Delta=h+Z$, the likelihood drop is $Z^2$
when $\Delta>0$ and $-h^2-2hZ$ when $\Delta\le0$.
\end{proof}



\section{Finite-sample illustration}
\label{sec:simulation}
 
\begin{figure}[htbp]
  \centering
  \includegraphics[width=0.8\textwidth]{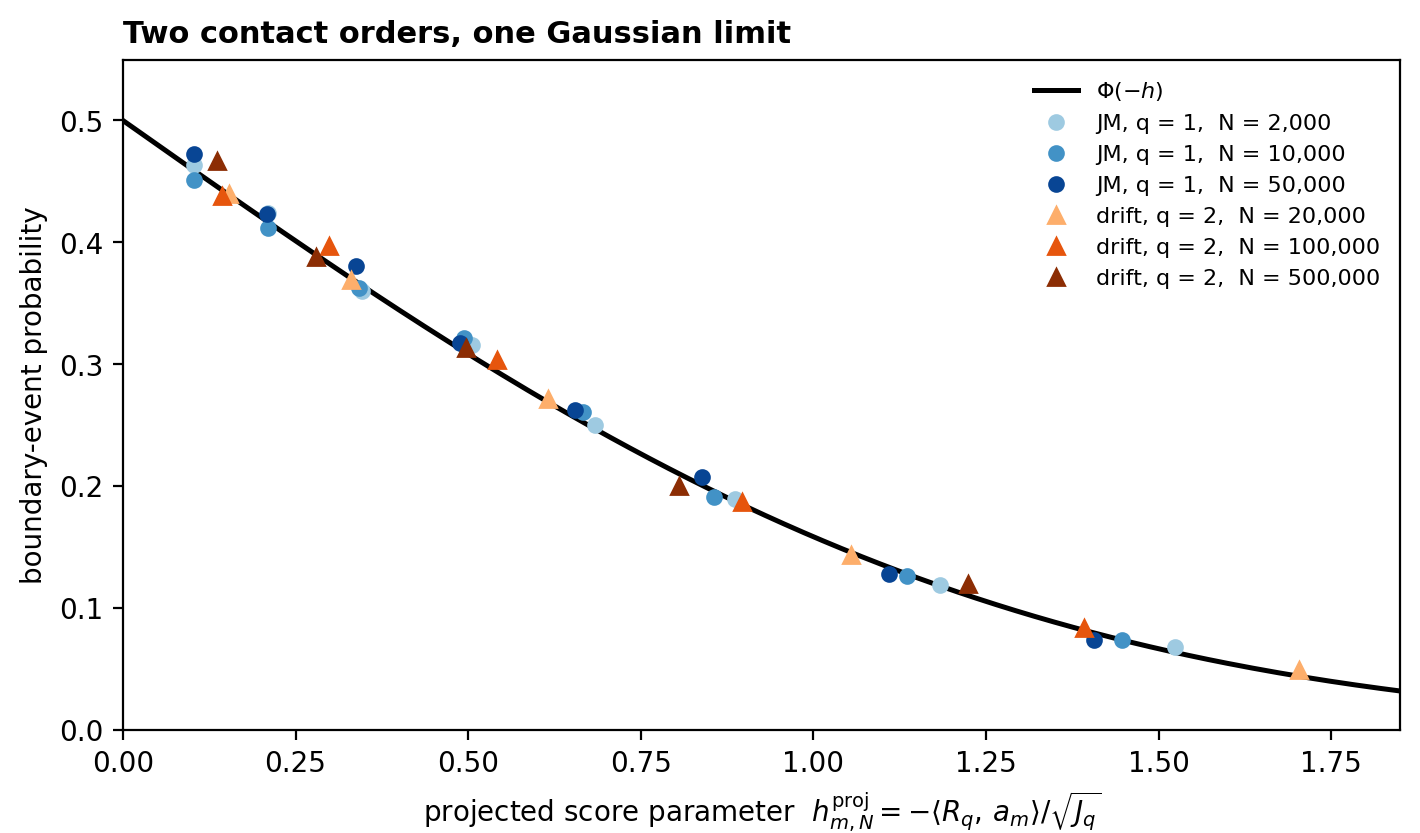}
  \caption{Master collapse across contact orders.  Simulated
  boundary-event probability against the finite-$N$ projected score parameter
  \eqref{eq:hexact} for the JM model ($q=1$, circles; the event is exactly
  $\{\widehat N=\infty\}$) and the effort-drift model ($q=2$, triangles;
  the score-sign boundary event of the limit experiment).  The black curve
  is $\Phi(-h)$.  JM: $N=2{,}000$, $10{,}000$, $50{,}000$ with
  $s\in\{0.5,0.8,1.1,1.4,1.7,2.0,2.4,2.8\}$ and $6{,}000$ replicates per
  point.  Effort drift: $N=20{,}000$, $100{,}000$, $500{,}000$ with
  $s\in\{1.6,2.1,2.6,3.1,3.6\}$ and $1{,}800$, $1{,}000$, $600$ replicates
  per point respectively.}
  \label{fig:mastercollapse}
\end{figure}

\begin{figure}[htbp]
  \centering
  \includegraphics[width=\textwidth]{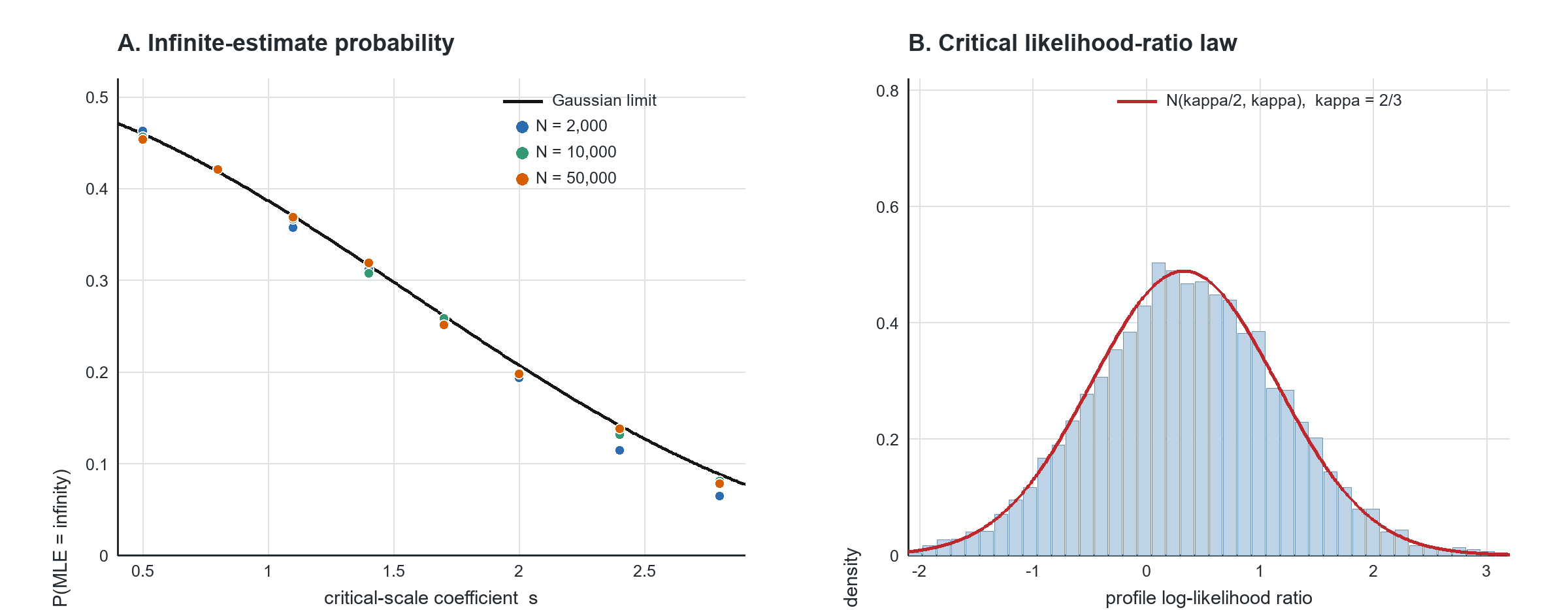}
  \caption{Critical-scale behavior in the JM model.  Left: probability of an
  infinite MLE, at $N=2{,}000$, $10{,}000$, $50{,}000$ and
  $s\in\{0.5,0.8,1.1,1.4,1.7,2.0,2.4,2.8\}$.  Right: distribution of the profile log-likelihood ratio comparing the
finite-population model with its limiting model at $s=2$,
  $N=50{,}000$.  Both panels use $6{,}000$ replicates per point.}
  \label{fig:simulation}
\end{figure}

\begin{figure}[htbp]
  \centering
  \includegraphics[width=0.92\textwidth]{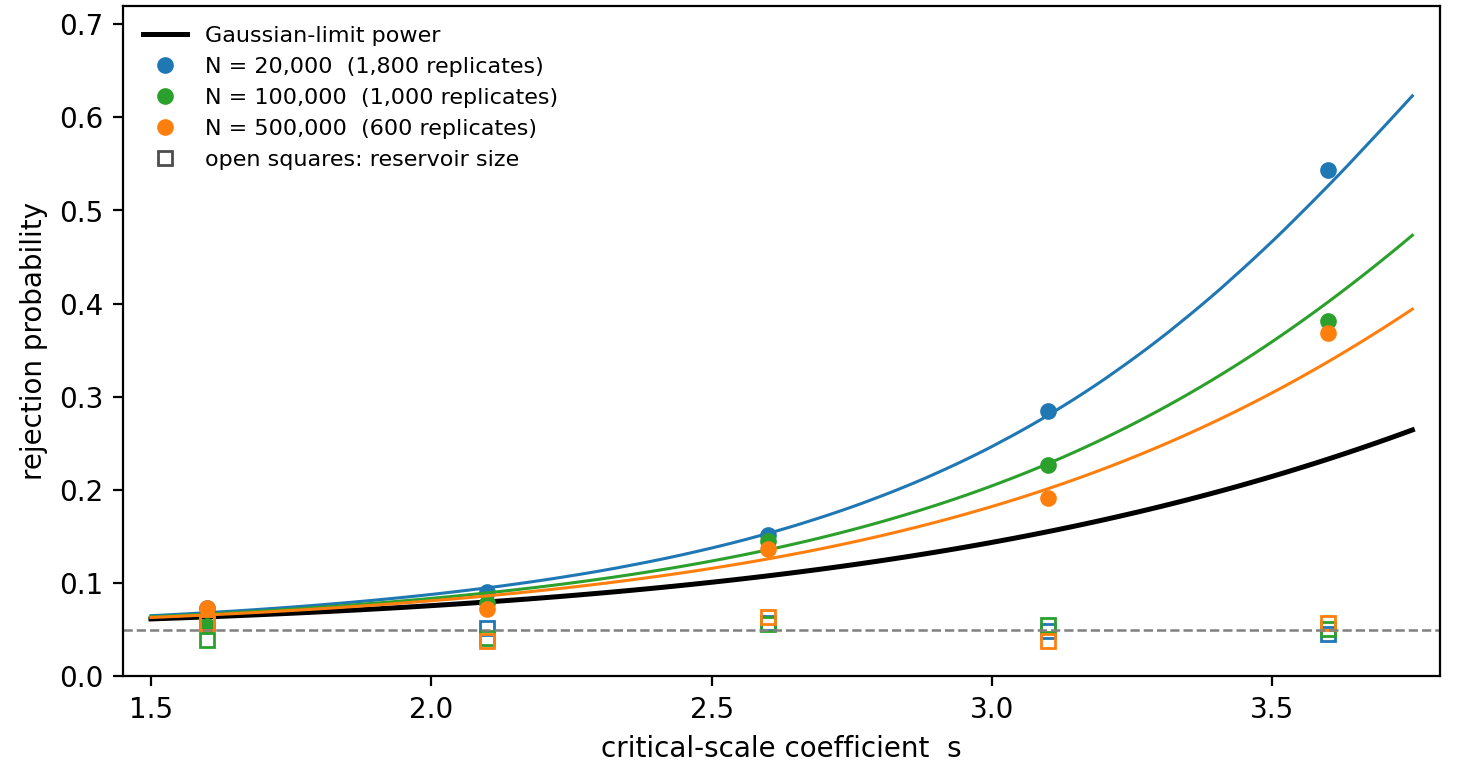}
  \caption{Size and power at the $N^{4/5}$ boundary in the effort-drift
  model.  Filled circles show rejection probabilities under matched finite
  alternatives, open squares show rejection probabilities under the
  limiting-model, the black curve is the critical Gaussian-limit power, and the
  colored curves are finite-$N$ predictions from the exact moments \eqref{eq:exactmoments}.  Designs are $N=20{,}000$, $100{,}000$,
  $500{,}000$ with $s\in\{1.6,2.1,2.6,3.1,3.6\}$, giving $m$ between
  $4{,}415$ and $130{,}460$, and $1{,}800$, $1{,}000$, $600$ replicates per
  point respectively.}
  \label{fig:driftsimulation}
\end{figure}

This section presents three simulation experiments: a joint check that both
contact orders are organized by a single information parameter
(Figure~\ref{fig:mastercollapse}), and separate checks of the
Jelinski-Moranda and effort-drift experiments
(Figures~\ref{fig:simulation} and~\ref{fig:driftsimulation} respectively).  
 
\paragraph{Simulation design.}
The simulations were conducted in the following manner, which
uses only elementary quantities defined above; no further calibration is
involved.
 
\emph{Jelinski-Moranda.}  Fix $N$ and $s$, and set
$m=\lfloor sN^{2/3}\rfloor$.  Draw $X_j\sim\Exp\{\phi(N-j)\}$
independently for $j=0,\ldots,m-1$; the rate $\phi$ may be set to one
without loss of generality, since every statistic used here is invariant
under a common rescaling of the $X_j$.  An infinite maximum-likelihood
estimate is recorded exactly through the
criterion of Proposition~\ref{prop:existence}: $\widehat N=\infty$ if and
only if $W_m=\sum_jjX_j/\sum_jX_j$ does not exceed $(m-1)/2$.  The profile
log-likelihood ratio is evaluated in closed form from the exact identity
\eqref{eq:exactprofile}.
 
\emph{Effort drift.}  Fix $N$ and $s$, and set $m=\lfloor sN^{4/5}\rfloor$
and $t_j=j/m$.  Under the limiting model, draw $X_j\sim\Exp(1)$, corresponding
to $\eta=\gamma=0$ in \eqref{eq:driftmodel}.  Under the matched finite
alternative of Theorem~\ref{thm:tangent}, draw $X_j\sim\Exp(e^{d_j})$,
where $d=(I_m-\Pi_m)a_m(1/N)$ is the residual of the log-rate profile on
$\operatorname{span}\{1,t\}$; equivalently $X_j=e^{-d_j}E_j$ with $E_j$
independent $\Exp(1)$.  Absorbing the affine component of $\log(1-j/N)$
into the unknown level and drift in this way leaves only the observable
residual log-rate, which is the point of the construction.  The test
statistic is the projected score
$\Delta_m^{\perp}=\sum_jR_{2,m,j}^{I}(\lambda_{j,0}X_j-1)/\sqrt{\mathcal
J_{2,m}}$ (with $\lambda_{j,0}$ denoting the rates under the limiting model), rejecting at one-sided
level $0.05$ when $\Delta_m^{\perp}>z_{0.95}$.
 
\emph{Reported quantities.}  Rejection probabilities and boundary-event
probabilities are Monte Carlo averages over independent replicates; the
replicate counts and design grids are given in the figure captions.
Standard errors are binomial, so a probability near $0.3$ estimated from
$600$ replicates carries a standard error of about $0.019$, which is the
resolution at which the figures should be read.

\paragraph{The two experiments share one limit.}
Figure~\ref{fig:mastercollapse} superposes both simulations on a single
axis.  Throughout, the horizontal coordinate is the finite-$N$ projected
score parameter
\begin{equation}
 h_{m,N}^{\rm proj}
 =-\frac{\ip{R_{q,m}^{I}}{a_m(1/N)}}{\sqrt{\mathcal J_{q,m}}},
 \qquad
 a_m(1/N)=\{\log(1-j/N):j<m\},
\label{eq:hexact}
\end{equation}
the standardized coefficient of the full finite-$N$ depletion profile in
the leading residual direction $R_{q,m}^{I}$.  Because $R_{q,m}^{I}$ is
orthogonal to $\operatorname{col}(Z_m)$, the projected and unprojected
profiles give the same inner product in \eqref{eq:hexact}.  In both models
$I(\theta)\equiv1$ on the log-rate scale, so $D_m=I_m$ and the Fisher and
Euclidean inner products coincide; for JM the residual is the centered
index, with $\mathcal J_{1,m}=V_m$, and for the effort-drift model it is
\eqref{eq:quadraticresidual}.

The
vertical coordinate is the simulated probability of the boundary event: for
JM this is $\{\Delta_m\le0\}$, which by
Proposition~\ref{prop:existence} is exactly $\{\widehat N=\infty\}$; for
the drift model it is $\{\Delta^{\perp}_m\le0\}$, the boundary event of the
limit experiment, in which the constrained maximizer is
$(\Delta^{\perp})_+$. 

All six design families (i.e., two contact orders, three population sizes
each) collapse on the single curve $\Phi(-h)$.

Two features of this collapse should be noted.  First,  the
choice of horizontal coordinate:  One could plot instead against the
limiting information parameter, which for the effort-drift model is
$s^{5/2}/\sqrt{720}$ and depends on $s$ alone; the two coordinates agree
when $m/N$ is small but diverge otherwise, because $h^{\rm proj}_{m,N}$
retains the cubic and higher terms of $\log(1-j/N)$, which project onto
$R^{I}_{2,m}$ with a common sign.  The most informative simulation in this
respect is $N=20{,}000$ with $s=3.6$, where $m/N\approx0.5$: there the
quadratic term alone contributes exactly the limiting value $0.92$, the
remaining terms carry $h^{\rm proj}_{m,N}$ to $1.70$, and the simulated
boundary probability follows the latter.  Plotted against the limiting
parameter, that point would miss the curve by a wide margin.  The
one-sided Gaussian description therefore continues to hold well beyond the
regime in which the analytic approximation to $h$ is accurate: what
degrades as $m/N$ grows is the approximation of the depletion profile by
its leading term, not the Gaussian structure built upon it.

The second feature is that the two experiments reach the same limit by
different routes.  The Jelinski-Moranda model has contact order $q=1$
against the nuisance space $\operatorname{span}\{1\}$, giving
$h_{m,N}\asymp m^{3/2}/N$ and the critical scale $N^{2/3}$; the
effort-drift model has effective contact order $q=2$ against
$\operatorname{span}\{1,j\}$, giving $h_{m,N}\asymp m^{5/2}/N^2$ and the
critical scale $N^{4/5}$.  The residuals, the exponents, and the required
sample sizes all differ (and at $N=500{,}000$ the two critical scales differ
by more than an order of magnitude) yet once each experiment is expressed
on its own $h$ scale the distinction disappears entirely, and both are the
one-sided Gaussian shift experiment $\mathcal N(h,1)$ with $h\ge0$.  The
contact order governs how a model approaches the boundary, not what is
found on arrival: the same boundary probability $\Phi(-h)$, the same power
envelope, and the same one-sided calibration apply to both.  This is the
sense in which contact order classifies rather than merely labels.

\paragraph{Jelinski-Moranda.}
Figure~\ref{fig:simulation} compares the exact JM experiment with the
critical limit.  The left panel shows the simulated probability of an
infinite MLE against the limit $\Phi(-s^{3/2}/\sqrt{12})$.  The right panel
compares the finite-model profile log-likelihood ratio at $s=2$ with its
$\mathcal N(\kappa/2,\kappa)$ limit, where $\kappa=s^3/12=2/3$; the
simulated mean and variance are $0.35$ and $0.71$, against $\kappa/2=0.33$
and $\kappa=0.67$.

The visible finite-sample departure occurs for large $s$ and small $N$.
This is expected because the neglected analytic corrections are controlled
by $m/N\asymp sN^{-1/3}$; for example, at $N=2{,}000$ and $s=2.4$ this
ratio is about $0.19$.  The primary collapse is organized by $h_{m,N}$,
while $m/N$ is the second parameter controlling the accuracy of the sparse
approximation.

\paragraph{Effort drift.}
Figure~\ref{fig:driftsimulation} gives a separate check of the genuinely
second-order example \eqref{eq:driftmodel}.  This experiment is reported
through size and power rather than through the probability of an infinite
estimate, because currently Proposition~\ref{prop:existence} has no counterpart here:
its proof relies on the unimodality of the JM profile likelihood, and with
the two-dimensional nuisance $(\eta,\gamma)$ neither an exact criterion nor
a closed-form profile is available.  Since
$h_{m,N}\to s^{5/2}/\sqrt{720}$, the limiting power is
\[
 \Phi\left\{\frac{s^{5/2}}{\sqrt{720}}-z_{0.95}\right\},
\]
shown as the black curve.  The colored curves use the finite-$N$
prediction obtained from the exact first two moments of the test statistic.
Under the matched alternative $\lambda_{j,0}X_j=e^{-d_j}E_j$ with $E_j$
independent $\Exp(1)$, and since $\sum_{j<m}R_{q,m,j}^{I}=0$,
\begin{equation}
 \Ee\Delta_m^{\perp}
 =\frac{\sum_{j<m}R_{q,m,j}^{I}\{e^{-d_j}-1\}}{\sqrt{\mathcal J_{q,m}}},
 \qquad
 \var\Delta_m^{\perp}
 =\frac{\sum_{j<m}(R_{q,m,j}^{I})^2e^{-2d_j}}{\mathcal J_{q,m}},
\label{eq:exactmoments}
\end{equation}
exactly, and the plotted curves are
$\Phi[\{\Ee\Delta^{\perp}_m-z_{0.95}\}/\{\var\Delta^{\perp}_m\}^{1/2}]$.
The parameter \eqref{eq:hexact} is the first-order term of the mean in
\eqref{eq:exactmoments}, since $e^{-d}-1=-d+d^2/2-\cdots$, and the variance
departs from one by $O(\max_j|d_j|)$; across the designs shown
$\max_j|d_j|\le0.045$, and the two predictions differ by at most $0.004$ in
rejection probability, an order of magnitude below the Monte Carlo error.

The colored predictions explain the slower approach to the black limit in
this example.  At the second-order boundary $m/N\asymp sN^{-1/5}$, rather
than $sN^{-1/3}$ as in JM, so the sparse approximation is reached more
slowly in $N$.  The simulations track the finite-$N$ prediction while the
latter approaches the $s^{5/2}/\sqrt{720}$ limit as $N$ grows.

The figure uses the oracle statistic, which depends on the unknown parameters
under the limiting model.  Substituting their maximum-likelihood estimates
under the limiting model, as in \eqref{eq:tangentfeasible}, leaves the
conclusions unchanged.  Evaluated on the same simulated samples, the feasible and
oracle statistics differ by a standard deviation of $0.024$ at the
smallest design considered and $0.004$ at the largest, against a statistic
of unit variance in both cases; the discrepancy decays at the rate
$m^{-1/2}$ predicted by Theorem~\ref{thm:tangent}, and the two tests
reach the same decision in more than $99\%$ of replicates.

\section{Discussion}

The $N^{2/3}$ law is not a peculiarity of an old software-reliability
model.  It is the $q=1$ member of an effective analytic-contact
classification, in which the squared norm of the first residual Taylor
trend determines both the information boundary and the local likelihood
experiment.  What matters is the nuisance structure, not the realism of the
model.  An unknown common level costs only the factor $q^2/(q+1)^2$ in
information and leaves the
exponent alone; profiling all polynomial trends below order $q$ costs the
larger but still constant fraction $\binom{2q}{q}^{-2}$.  Enlarging the
nuisance space, by contrast, can absorb a complete Taylor coefficient and
so change the exponent itself, as the effort-drift model does in moving the
boundary from $N^{2/3}$ to $N^{4/5}$.  Contact order is therefore a joint
property of the depletion profile and the nuisance space, and the
residual-norm calculus is what makes it computable.

The results also clarify what an infinite maximum-likelihood estimate does and does not mean.  It is not evidence that the population is literally infinite.  Near the boundary it is a typical  outcome with calculable probability $\Phi(-h)$, and under the exact limiting model it occurs one half of
the time, for every sample size.  Replacing infinity by a numerical cap or
by an unreported prior does not create information; it merely relocates the
difficulty.  Honest inference should report the one-sided likelihood
geometry directly, with two immediate consequences in practice.  The null
distribution of the likelihood-ratio statistic is
$\tfrac12\delta_0+\tfrac12\chi^2_1$,  and the
confidence set of Corollary~\ref{cor:confidence} is allowed to be unbounded,
being unbounded exactly when the data does not exclude the limiting model, with coverage $1-\alpha$ once $h\ge z_{1-\alpha/2}$ and
$1-\alpha/2$ at the boundary itself.  Such a set reports what the
experiment can support rather than what a chosen upper limit permits.

The information constants also carry a design implication.  Because
profiling all polynomial trends below order $q$ retains only
$\binom{2q}{q}^{-2}$ of the available information, plummeting  from $\tfrac1{4}$ at $q=1$
to $\tfrac1{400}$ at $q=3$, the cost of treating a covariate as an
unknown nuisance trend rather than measuring it grows very quickly.  An
analyst who suspects an effort drift and cannot exclude curvature in it is
already at $q=3$; recording testing hours, exposure, or sampling intensity
directly and entering them as known offsets recovers that factor in full.
When the depletion signal is weak, instrumenting the experiment is worth
considerably more than refining the analysis of an uninstrumented one.

The approach presented here
suggests a broader analytic theory of likelihoods under partial observation.
For example,  an invariance description for the
exponent: thinning, aggregation, interior censoring, and marginalization
each act on the residual analytic contrast, and the organizing question is
which of these operations preserve the contact order, which raise it, and
which lower it.  Thinning and aggregation appear to preserve it and to act
only on the constant, whereas marginalization to a statistic which is insensitive
at first order should raise it. A classification of this kind
would allow the exponent to be read off from the model, the nuisance
parameters, and the manner in which the data was degraded.  

\appendix

\section{Proof of Theorem~\ref{thm:lan}}
\label{app:lan}

Write $I_0=I(\eta_0)$, $g_j=g_{j,q,m}$, and
$\mathcal V_m=V_{q,m}$.  Fix $H<\infty$.  Uniformly for $0\le h\le H$,
put $\rho=\rho_m(h)$ and
$\widetilde a_j=\widetilde a_{j,m}(\rho)$.  Equation \eqref{eq:Vasym} gives
\[
 \rho
 =
 O\{m^{-1-1/(2q)}\},
 \qquad
 m\rho=O\{m^{-1/(2q)}\}\to0.
\]
Define the centered analytic remainder
\[
 \varepsilon_{j,m}
 :=
 \widetilde a_j+c\rho^q g_j,
 \qquad\text{so that}\qquad
 \widetilde a_j
 =
 -\frac{h}{\sqrt{I_0\mathcal V_m}}g_j+\varepsilon_{j,m},
 \qquad
 \sum_{j<m}\varepsilon_{j,m}=0,
\]
the zero sum holding because both $\widetilde a_j$ and $g_j$ are
centered.  The remainder is negligible in quadratic mean because Taylor's
bound on $a$ gives
\[
 \sup_{0\le h\le H}\sum_{j<m}\varepsilon_{j,m}^2
 =
 O\{\rho^{2q+2}m^{2q+3}\}
 =
 O(m^{-1/q})\to0.
\]
Consequently,
\[
 \max_{j<m}|\widetilde a_j|=O(m^{-1/2}),\qquad
 \sum_{j<m}\widetilde a_j=0,
 \qquad
 I_0\sum_{j<m}\widetilde a_j^2=h^2+o(1).
\]

Taylor's theorem gives
\[
\sum_{j<m}
 \{\ell_{\eta_0+\widetilde a_j}(Y_j)-\ell_{\eta_0}(Y_j)\}
 =
 \sum_{j<m}\widetilde a_j\dot\ell_{\eta_0}(Y_j)
 +
 \frac12\sum_{j<m}\widetilde a_j^2\ddot\ell_{\eta_0}(Y_j)
 +
 o_{\Pp}(1).
\]
Indeed, (S2) and the law of large numbers bound the summed third-order
remainder, uniformly in $h\le H$, by
\[
 O_{\Pp}\{\max_j|\widetilde a_j|\sum_j\widetilde a_j^2\}
 =O_{\Pp}(m^{-1/2}).
\]
A weighted law of large numbers and
$\Ee_{\eta_0}\ddot\ell_{\eta_0}=-I_0$ give
\[
 \frac12\sum_{j<m}\widetilde a_j^2\ddot\ell_{\eta_0}(Y_j)
 =
 -\frac{h^2}{2}+o_{\Pp}(1).
\]
The score contribution of $\varepsilon_{j,m}$ has variance
$I_0\sum_j\varepsilon_{j,m}^2=o(1)$, and hence the linear term is
$h\Delta_m+o_{\Pp}(1)$.

To profile the common parameter, reparametrize it as
$\eta+\bar a_m(\rho)$ and allow a local displacement
$b/\sqrt{mI_0}$.  The nuisance and trend central sequences are
\[
 \Delta_m^{\parallel}
 =\frac1{\sqrt{mI_0}}\sum_{j<m}\dot\ell_{\eta_0}(Y_j),
 \qquad
 \Delta_m=
 -\frac1{\sqrt{I_0\mathcal V_m}}\sum_{j<m}g_j\dot\ell_{\eta_0}(Y_j),
\]
the superscript recording that the first points along the nuisance tangent
direction while the second is orthogonal to it.  The same Taylor
calculation, now uniform for bounded $b$ and $h$, gives
\[
 b\Delta_m^{\parallel}+h\Delta_m-\frac12(b^2+h^2)+o_{\Pp}(1).
\]
The expected cross information vanishes because
$\sum_j\widetilde a_j=0$.  The
observed cross term is
\[
 \frac{b}{\sqrt{mI_0}}\sum_j \widetilde a_j
 \{\ddot\ell_{\eta_0}(Y_j)+I_0\}=o_{\Pp}(1)
\]
uniformly for bounded $b,h$: its variance is
$O(m^{-1}\sum_j\widetilde a_j^2)=O(m^{-1})$.
Assumption (S3), followed by the score expansion, makes the maximizing
$b$ tight.  Maximizing over $b$ under both the finite and limiting-model models
cancels $b\Delta_m^{\parallel}-b^2/2$ and proves
\eqref{eq:profilelan}.

For the null limit, the largest deterministic weight in $\Delta_m$ is
\[
 \max_{j<m}\frac{|g_j|}{\sqrt{\mathcal V_m}}=O(m^{-1/2}).
\]
Assumption (S2) therefore verifies the Lindeberg condition, and the
weighted central limit theorem gives
$\Delta_m\Rightarrow\mathcal N(0,1)$.

It remains to calculate the distribution under
\eqref{eq:matchedalternative}.  Differentiating the density under the
integral, justified by (S1)-(S2), yields uniformly for small $u$,
\[
 \Ee_{\eta_0+u}\dot\ell_{\eta_0}(Y)
 =
 I_0u+O(u^2),
 \qquad
 \Var_{\eta_0+u}\{\dot\ell_{\eta_0}(Y)\}
 =
 I_0+o(1).
\]
Substituting the centered shifts $\widetilde a_j$ shows directly that
\[
 \Ee(\Delta_m)=h_0+o(1),
 \qquad
 \Var(\Delta_m)=1+o(1).
\]
The same maximal-weight bound and the uniform $(2+\delta)$ moment in (S2)
give the Lindeberg condition under the alternative.  The centered
triangular-array central limit theorem proves
\eqref{eq:alternativeclt}.

The likelihood expansion under the alternative follows by the same Taylor
calculation: (S2) is uniform over the $O(m^{-1/2})$ parameter shifts, the
quadratic term remains $-h^2/2+o_{\Pp}(1)$, and the linear term has the
alternative central-sequence limit just derived.  This also establishes the
compact-uniform statement without an appeal to a change-of-measure lemma.

Finally, let $\widehat\eta_0$ be the homogeneous MLE.  By (S3),
$\widehat\eta_0-\eta_0=O_{\Pp}(m^{-1/2})$.  Taylor expansion of the score in
\eqref{eq:feasiblescore} gives, using $\sum_jg_j=0$,
\[
 \sum_jg_j\{\dot\ell_{\widehat\eta_0}(Y_j)
 -\dot\ell_{\eta_0}(Y_j)\}
 = (\widehat\eta_0-\eta_0)
 \sum_jg_j\{\ddot\ell_{\eta_0}(Y_j)+I_0\}
 +o_{\Pp}(\sqrt{\mathcal V_m}).
\]
The centered sum on the right is $O_{\Pp}(\sqrt{\mathcal V_m})$, so after
normalization the first term is $o_{\Pp}(1)$; the second-order score
remainder is also $O_{\Pp}(m^{-1/2})$ after normalization.  Continuity of
$I$ completes \eqref{eq:plugin}.  The identical argument applies under
matched local alternatives.
\qed

\section{Proof of Theorem~\ref{thm:tangent}}
\label{app:tangent}

Write $R=R_{q,m}^{I}$, $\mathcal J=\mathcal J_{q,m}$,
$I_j=I(\theta_{j,0})$, and $D=D_m$.  The weighted Gram condition and the
uniform bounds on $I_j$ imply that $\Pi_{m,I}$ has uniformly bounded
Euclidean operator norm.  Since every active Taylor trend below order $q$
belongs to $\operatorname{col}(Z_m)$, projection of the analytic expansion
gives, uniformly for bounded $h$,
\[
 d_m\{\rho_m(h)\}
 =-c\rho_m(h)^qR+\varepsilon_m
 =-\frac{h}{\sqrt{\mathcal J}}R+\varepsilon_m,
\]
where
\[
 \varepsilon_m^\mathsf TD\varepsilon_m
 =O\{\rho_m(h)^{2q+2}m^{2q+3}\}
 =O(m^{-1/q})=o(1).
\]
Condition \eqref{eq:effectiveconditions} also gives
\[
 \max_{j<m}|d_{m,j}|=O(m^{-1/2}),
 \qquad
 d_m^\mathsf TD_md_m=h^2+o(1),
 \qquad
 Z_m^\mathsf TD_md_m=0.
\]
Taylor expansion around the nonidentical baseline parameters $\theta_{j,0}$
therefore yields
\[
 \sum_{j<m}\{\ell_{\theta_{j,0}+d_{m,j}}(Y_j)
                 -\ell_{\theta_{j,0}}(Y_j)\}
 =h\Delta_m^\perp-\frac{h^2}{2}+o_{\Pp}(1).
\]
Indeed, the score contribution of $\varepsilon_m$ has variance $o(1)$, a weighted law
of large numbers replaces $\ddot\ell_{\theta_{j,0}}$ by $-I_j$, and (S2)
bounds the summed third-order remainder by
$O_{\Pp}(\max_j|d_{m,j}|\sum_jd_{m,j}^2)=o_{\Pp}(1)$.

It remains to check profiling.  After absorbing
$\Pi_{m,I}a_m(\rho)=Z_m\gamma_m(\rho)$ into $\beta$, let the nuisance
parameter move from $\beta_{0,m}$ by $b/\sqrt m$, with bounded
$b\in\R^p$.  The deterministic cross information is zero because
$Z_m^\mathsf TD_md_m=0$.  The observed cross term is
\[
 \frac1{\sqrt m}\sum_{j<m}(z_{j,m}^\mathsf Tb)d_{m,j}
 \{\ddot\ell_{\theta_{j,0}}(Y_j)+I_j\},
\]
whose variance is $O(m^{-1}\sum_jd_{m,j}^2)=O(m^{-1})$.  Thus it is
$o_{\Pp}(1)$ uniformly for bounded $b,h$.  The nuisance quadratic is the
same under the limiting-model and finite profiles and cancels on maximization.
Root-$m$ localization of the MLE under the limiting model makes the maximizing $b$ tight,
which proves \eqref{eq:tangentlan}.

The maximal-weight condition in \eqref{eq:effectiveconditions}, together
with (S2), verifies Lindeberg and yields the null CLT.  Under the matched
alternative, uniformly over the baseline curve,
\[
 \Ee_{\theta+u}\dot\ell_\theta=I(\theta)u+O(u^2),
 \qquad
 \Var_{\theta+u}(\dot\ell_\theta)=I(\theta)+o(1).
\]
Substitution of $d_m=-h_0R/\sqrt{\mathcal J}+\varepsilon_m$ gives mean
$h_0+o(1)$ and variance $1+o(1)$ for $\Delta_m^\perp$; the same weighted
CLT applies.  All Taylor remainders remain compact-uniform under these local
alternatives.

For feasibility, continuous differentiability of $I$, bounded regressors,
and root-$m$ consistency give
\[
 \|\widehat D_m-D_m\|_{\rm op}=O_{\Pp}(m^{-1/2}),
 \qquad
 \|\widehat R_{q,m}^{I}-R_{q,m}^{I}\|
   =O_{\Pp}(m^q),
 \qquad
 \frac{\widehat{\mathcal J}_{q,m}}{\mathcal J_{q,m}}
   =1+O_{\Pp}(m^{-1/2}).
\]
The second bound is smaller by $m^{-1/2}$ than
$\|R_{q,m}^{I}\|=O(m^{q+1/2})$, so replacing the oracle residual and
normalizer changes the standardized score by $o_{\Pp}(1)$.  Expanding
$\dot\ell_{\widehat\theta_{j,0}}$ around $\theta_{j,0}$, the deterministic
first-order nuisance term vanishes under the estimated Fisher metric because
$Z_m^\mathsf T\widehat D_m\widehat R_{q,m}^{I}=0$.  The difference between
$\widehat D_m$ and $D_m$, the centered observed-Hessian term, and the
second-order score remainder are all $O_{\Pp}(m^{-1/2})$ after division by
$\sqrt{\mathcal J_{q,m}}$.  This proves \eqref{eq:tangentfeasible}.

Finally, the same quadratic expansion gives
\[
 \KL\!\left(
  \bigotimes_{j<m}P_{\theta_{j,0}+d_{m,j}(\rho)},
  \bigotimes_{j<m}P_{\theta_{j,0}}
 \right)
 =\frac{h_m^2}{2}\{1+o(1)\}.
\]
Pinsker's inequality proves merging when $h_m\to0$; the Gaussian limit gives
the critical experiment, and a centered-score threshold proves consistent
separation when $h_m\to\infty$ while $m\rho\to0$.
\qed

\section{Proof of Corollary~\ref{cor:trichotomy}}

Take the closest homogeneous parameter
\[
 \eta_m^\star
 =
 \eta+\frac1m\sum_{j<m}a(\rho j)
 =
 \eta-\frac{c\rho^q}m\sum_{j<m}j^q+o(\rho^qm^q).
\]
The residual parameter displacement has squared information norm
$h_m^2\{1+o(1)\}$.  The second-order log-density expansion in (S1), with
the uniform integrability supplied by (S2), gives
\[
 \KL\!\left(
   \bigotimes_{j<m}P_{\eta+a(\rho j)},
   P_{\eta_m^\star}^{\otimes m}
 \right)
 =
 \frac{h_m^2}{2}\{1+o(1)\}.
\]
If $h_m\to0$, Pinsker's inequality proves total-variation merging.
The critical assertion is Theorem~\ref{thm:lan}.  If $h_m\to\infty$ while
$m\rho\to0$, the centered score has null mean zero, asymptotic variance one,
and finite-model mean $h_m\{1+o(1)\}$.  A threshold diverging more slowly than
$h_m$ produces both type-I and type-II errors tending to zero.
\qed

\section{Proof of Proposition~\ref{prop:exactexp} and
Corollary~\ref{cor:lrlimit}}

For the finite model, profiling $\alpha$ gives
\[
 \sup_\alpha\ell_f(\alpha)
 =
 C_m+\sum_{j<m}\log f_j
 -m\log\left(\sum_{j<m}f_jX_j\right),
\]
whereas profiling the common rate \(\beta\) under the limiting model gives
\[
 \sup_\beta\ell_\infty(\beta)
 =
 C_m-m\log\left(\sum_{j<m}X_j\right).
\]
Substitute $X_j=b_{j,N}E_j/\alpha$, use
$b_{j,N}=\bar b_{m,N}\sigma_{j,N}$, and cancel the terms involving
$\bar b_{m,N}$ to obtain
\eqref{eq:exactprofile}.

Both limits now follow from a single expansion, applied with reciprocal
weightings.

\begin{lemma}[Profiled exponential contrast]
\label{lem:contrast}
Let $c_0,\ldots,c_{m-1}$ be positive constants with $m^{-1}\sum_jc_j=1$,
put $\delta_j=1-c_j$, and let $G_0,\ldots,G_{m-1}$ be independent
$\Exp(1)$.  If $\max_j|\delta_j|\to0$ and
$\sum_j\delta_j^2\to\kappa\in[0,\infty)$, then
\[
 \Lambda(c)
 :=
 -\sum_{j<m}\log c_j
 +m\log\frac{\sum_jc_jG_j}{\sum_jG_j}
 =
 \frac12\sum_j\delta_j^2
 -\sum_j\delta_j(G_j-1)
 +o_{\Pp}(1)
 \;\Longrightarrow\;
 \mathcal N(\kappa/2,\kappa).
\]
\end{lemma}

\begin{proof}
Since $\sum_j\delta_j=0$ and
$\sum_j|\delta_j|^3\le\max_j|\delta_j|\sum_j\delta_j^2\to0$,
\[
 -\sum_j\log c_j
 =
 -\sum_j\log(1-\delta_j)
 =
 \sum_j\delta_j+\frac12\sum_j\delta_j^2
 +O\Bigl(\sum_j|\delta_j|^3\Bigr)
 =
 \frac12\sum_j\delta_j^2+o(1).
\]
Next, $\sum_jc_jG_j/\sum_jG_j=1-S_m/\sum_jG_j$ with
$S_m=\sum_j\delta_j(G_j-1)=O_{\Pp}(\kappa^{1/2})$, again using
$\sum_j\delta_j=0$.  Since $m^{-1}\sum_jG_j\to1$ in probability,
\[
 m\log\Bigl(1-\frac{S_m}{\sum_jG_j}\Bigr)
 =
 -\frac{mS_m}{\sum_jG_j}
 +m\cdot O_{\Pp}\Bigl(\frac{S_m^2}{m^2}\Bigr)
 =
 -S_m+o_{\Pp}(1).
\]
The weights $\delta_j$ satisfy the Lindeberg condition because
$\max_j|\delta_j|\to0$, so $S_m\Rightarrow\mathcal N(0,\kappa)$, proving
the displayed limit.
\end{proof}

Under the finite model, \eqref{eq:exactprofile} says exactly
$\Lambda_{m,N}=\Lambda(\sigma)$ with $c_j=\sigma_{j,N}$ and $G_j=E_j$;
since $1-\sigma_{j,N}=-c\rho^qg_{j,q,m}\{1+o(1)\}$ elementwise,
$\sum_j(1-\sigma_{j,N})^2=\kappa_{m,N}\{1+o(1)\}$, and
Lemma~\ref{lem:contrast} gives \eqref{eq:lrfinite}.  Under the limiting model,
with $\tau_{j,N}=f(j/N)/\{m^{-1}\sum_if(i/N)\}$ and $G_j=\beta X_j$, the
same profiled likelihood is exactly $\Lambda_{m,N}=-\Lambda(\tau)$, and
$\sum_j(1-\tau_{j,N})^2=\kappa_{m,N}\{1+o(1)\}$, giving
\eqref{eq:lrnull}.  The $\pm\kappa/2$ symmetry of
Corollary~\ref{cor:lrlimit} is thus structural: the two hypotheses apply
the same contrast with reciprocal weightings.
\qed

\section{Proof of Theorem~\ref{thm:jm} and
Corollary~\ref{cor:mlelaw}}
The JM model is \eqref{eq:exparray} with
$\alpha=\phi N$, $f(x)=1-x$, and $\rho=N^{-1}$.  Theorem~\ref{thm:lan}
therefore gives the local likelihood process and
\[
 \Delta_m\Rightarrow\mathcal N(h,1).
\]
The scale-invariant form \eqref{eq:Deltajm} differs from the efficient central sequence by $o_{\Pp}(1)$ because $m/S_m$ consistently estimates the matched rate under the limiting model.  The merging and separation statements follow from
Corollary~\ref{cor:trichotomy}.

The maximum-likelihood result requires one further property of the exact
profile likelihood, which we establish directly.

\begin{lemma}[Unimodality of the JM profile likelihood]
\label{lem:unimodal}
On $\rho\in[0,1/m]$ the profile log-likelihood \eqref{eq:jmprofile} has at
most one stationary point, and every stationary point is a strict local
maximum.  It is therefore unimodal, with a unique maximizer.
\end{lemma}

\begin{proof}
Write
\[
 \varphi(\rho)
 =
 \frac{d}{d\rho}\pl_m(\rho)
 =
 m\frac{W_m}{1-\rho W_m}
 -
 \sum_{j=0}^{m-1}\frac{j}{1-\rho j}
\]
for the profile score obtained in the proof of
Proposition~\ref{prop:existence}, and put $v_j=j/(1-\rho j)$ for $j<m$ and
$\bar v=W_m/(1-\rho W_m)$; these are finite because $\rho\le1/m$ and
$W_m\le m-1$.  Differentiating once more,
\[
 \varphi'(\rho)
 =
 \frac{mW_m^2}{(1-\rho W_m)^2}
 -\sum_{j=0}^{m-1}\frac{j^2}{(1-\rho j)^2}
 =
 m\bar v^2-\sum_{j=0}^{m-1}v_j^2 .
\]
At a stationary point $\varphi(\rho)=0$, that is $\sum_{j<m}v_j=m\bar v$,
so $\bar v$ is the average of the $v_j$ and
\[
 \varphi'(\rho)=-\sum_{j=0}^{m-1}(v_j-\bar v)^2<0,
\]
strictly, since $j\mapsto j/(1-\rho j)$ is strictly increasing and the
indices are distinct for $m\ge2$.  Every stationary point is therefore a
strict local maximum, and a differentiable function on an interval cannot
have two strict local maxima without an intervening stationary point that
is not one.  Hence there is at most one.
\end{proof}

For the MLE result, map the integer parameter set to
\[
 \mathcal H_m
 =
 \left\{\frac{\sqrt{V_m}}{n}:n=m,m+1,\ldots\right\}\cup\{0\}.
\]
On every compact subset of $(0,\infty)$ the maximal spacing of
$\mathcal H_m$ tends to zero.  Lemma~\ref{lem:unimodal}, together with the
uniform compact likelihood expansion, makes
the maximizer tight in the $h$ scale: the profile likelihood has a single
mode, so it cannot exceed near any larger $u$ the value it attains near
$\Delta_m$.  The argmax continuous mapping theorem
then gives
\[
 \widehat h_m
 \Rightarrow
 \mathop{\rm argmax}_{u\ge0}
 \left\{u\Delta-\frac{u^2}{2}\right\}
 =
 (\Delta)_+.
\]
The transformation in \eqref{eq:Nhat} follows, with the convention that
$1/0=\infty$.
\qed

\bibliographystyle{plainnat}
\small
\setlength{\bibsep}{0pt}

\end{document}